%% file: bf_T.tex
\documentclass[11pt, DIV10,a4paper]{article}
\usepackage{natbib}
\newcommand {\ctn}{\citet} 
\usepackage{graphicx,subfigure,amsmath,latexsym,amssymb}
\usepackage{float,epsfig,multirow,rotating,times}
\usepackage{upgreek,wrapfig}
\usepackage{comment}

\newcommand{\btheta}{\boldsymbol{\theta}}

\newcommand{\bbeta}{\boldsymbol{\beta}}

\newcommand{\bxi}{\boldsymbol{\xi}}
\newcommand{\bGamma}{\boldsymbol{\Gamma}}
\newcommand{\bTheta}{\boldsymbol{\Theta}}

\newcommand{\bC}{\boldsymbol{C}}

\newcommand{\bI}{\boldsymbol{I}}

\newcommand{\bZ}{\boldsymbol{Z}}
\newcommand{\bz}{\boldsymbol{z}}

\newtheorem{theorem}{Theorem}

\newtheorem{corollary}[theorem]{Corollary}

\newtheorem{lemma}[theorem]{Lemma}

\newtheorem{remark}[theorem]{Remark}

\newenvironment{proof}[1][Proof]{\textbf{#1.} }{\ \rule{0.5em}{0.5em}}

\newcommand{\boldcal}[1]{\mbox{\boldmath{$\mathcal #1 $}}}

\numberwithin{equation}{section}
\numberwithin{algo}{section}
\numberwithin{table}{section}
\numberwithin{figure}{section}

\usepackage{setspace}
\usepackage[mathscr]{euscript}
\usepackage[margin=1in]{geometry}
\begin{document}
\normalsize

\title{\vspace{-0.8in}
Asymptotic Theory of Bayes Factor in Stochastic Differential Equations: Part II}
\author{Trisha Maitra and Sourabh Bhattacharya\thanks{
Trisha Maitra is a PhD student and Sourabh Bhattacharya 
is an Associate Professor in
Interdisciplinary Statistical Research Unit, Indian Statistical
Institute, 203, B. T. Road, Kolkata 700108.
Corresponding e-mail: sourabh@isical.ac.in.}}
\date{\vspace{-0.5in}}
\maketitle%

\begin{abstract}

The problem of model selection in the context of a system of stochastic differential equations ($SDE$'s) 
has not been touched upon in the literature. Indeed, properties of Bayes factors have not been studied even 
in single $SDE$ based model comparison problems.

In this article, we first develop an asymptotic theory of Bayes factors when two $SDE$'s are compared, 
assuming the time domain expands. Using this we then develop an asymptotic theory of Bayes factors 
when systems of $SDE$'s are compared, 
assuming that the number of equations in each system, as
well as the time domain, increase indefinitely. 
Our asymptotic theory covers situations when the observed processes associated with the
$SDE$'s are independently and identically distributed ($iid$), as well as when they are independently 
but not identically distributed (non-$iid$). Quite importantly, we allow inclusion of available time-dependent 
covariate information into each $SDE$ through a multiplicative factor of the drift function in a random effects set-up; 
different initial values for the $SDE$'s are also permitted.

Thus, our general model-selection framework includes simultaneously the variable selection problem associated with
time-varying covariates, as well as choice of the part of the drift function free of covariates.
It is to be noted that given that the underlying process is wholly observed, the diffusion coefficient becomes known,
and hence is not involved in the model selection problem.

For both $iid$ and non-$iid$ set-ups we establish almost sure exponential convergence
of the Bayes factor. As we show, the Bayes factor is inconsistent for comparing individual $SDE$'s,
in the sense that the log-Bayes factor converges only in expectation, while
the relevant variance does not converge to zero.
Nevertheless, 
it has been possible to exploit this result to establish
almost sure exponential convergence of the Bayes factor when, in addition, the number of individuals 
are also allowed to increase
indefinitely.

We carry out simulated and real data analyses to demonstrate that Bayes factor is a suitable candidate
for covariate selection in our $SDE$ models even in non-asymptotic situations.
\\[2mm]
 
\end{abstract}

\section{Introduction}
\label{sec:intro}

Stochastic differential equations ($SDE$'s) have important standing in statistical applications 
where ``within" subject variability is caused by some random component
varying continuously in time. It also seems worthwhile to incorporate available time-dependent covariate
information into the subject-wise $SDE$'s.
Apart from the covariates there may also be random effects
associated with the individuals, which may be useful in modeling variabilities between the individuals.

$SDE$-based models with time-dependent covariates are considered in
\ctn{Zita11}, \ctn{Overgaard05}, \ctn{Leander15}; moreover, \ctn{Zita11} analyse their covariate-based $SDE$ model in 
the hierarchical Bayesian paradigm. In the literature, random effects $SDE$ models without covariates seem to be 
more popular than those based on covariates. 
A brief overview of random effects $SDE$ models is provided in
\ctn{Maud12} who undertake theoretical and classical asymptotic 
investigation of a class of
random effects models based on $SDE$'s. Specifically, they model the $i$-th individual
by
\begin{equation}
d X_i(t)=b(X_i(t),\phi_i)dt+\sigma(X_i(t))dW_i(t),
\label{eq:sde_basic1}
\end{equation}
where, for $i=1,\ldots,n$, $X_i(0)=x^i$ is the initial value of the stochastic process $X_i(t)$, which 
is assumed to be continuously observed on the time interval $[0,T_i]$; $T_i>0$ assumed to be known.
The function $b(x,\varphi)$, which is the drift function, is a known, real-valued function on $\mathbb R\times\mathbb R^d$ 
($\mathbb R$ is the real line and $d$ is the dimension), and the function $\sigma:\mathbb R\mapsto\mathbb R$ is the known 
diffusion coefficient.
The $SDE$'s given by (\ref{eq:sde_basic1}) are driven by independent standard Wiener processes $\{W_i(\cdot);~i=1,\ldots,n\}$, 
and $\{\phi_i;~i=1,\ldots,n\}$, which are to be interpreted as the random effect parameters associated
with the $n$ individuals, which are assumed by \ctn{Maud12} to be independent of the Brownian motions and
independently and identically distributed ($iid$) random variables with some common distribution.
For the sake of convenience 
\ctn{Maud12} (see also \ctn{Maitra14a} and \ctn{Maitra14b}) assume $b(x,\phi_i)=\phi_ib(x)$. 
Thus, the random effect is a multiplicative factor of the drift function. In this work,
we generalize this to a random effects $SDE$ set-up consisting of time-dependent covariates. 

Note that model selection constitutes an important part of research in both Bayesian and 
classical paradigms; see, for example, 
\ctn{Dey00}, \ctn{Jiang07}, \ctn{Claeskens08}, \ctn{Muller13}.
In the case of $SDE$-based mixed effects models as well, model selection constitutes an important issue
involving the choice of the drift function and selection of the appropriate subset of (time-dependent) covariates.
Here Bayes factors are expected to play the central role as their effectiveness in model selection
in complex problems is well-established (see, for example, \ctn{Kass95} for a good account of Bayes factors).
Unavailability of closed form expressions in the traditional $SDE$ set-ups usually
prompt usage of numerical approximations based on Markov chain Monte Carlo or 
related criteria such as the Akaike Information Criterion (\ctn{Akaike73}) 
and Bayes Information Criterion (\ctn{Schwarz78}). For details, see, for example, \ctn{Fuchs13}, \ctn{Iacus08}.
But we are not aware of any research existing in the literature that attempts to address covariate selection in $SDE$'s.

We are also not aware of any existing literature on asymptotic investigation of Bayes factors in the $SDE$ context 
although \ctn{Siva02} present some asymptotic investigation of intrinsic
and fractional Bayes factors in the context of three specific diffusion models.
The only investigation available in this context seems to be 
that of \ctn{Maitra15a}, who model a multiplicative part of the drift function using time-varying covariates, and
address Bayes factor asymptotics in a general set-up consisting of the 
covariate selection problem as well as selection of 
the part of the drift function independent of the covariates. Different initial values
and domains of observations pertaining to different individuals, are also considered in their set-up.
Assuming that only the number of individuals increase without bound, \ctn{Maitra15a} establish
almost sure exponential convergence of Bayes factor in both $iid$ and non-$iid$ situations.
Here we recall that the $iid$ set-up is the case when there is no covariate associated with 
the model and when the initial values and the domains of observations
are the same for every individual. 
The non-$iid$ set-up, on the other hand, consists of time-varying covariates, different  
initial values and domains of observations; in this work we also consider random effects. 
Thus, unlike the $iid$ case, here the model selection problem also deals with covariate selection 
apart from selection of the part of the drift functions free of the covariates.

In this article, we prove almost sure exponential convergence of the relevant Bayes 
factors in both $iid$ and non-$iid$ cases, assuming that the number
of individuals, as well as the domains of observations, increase without bound.
Hence, for our current purpose, the asymptotic theory developed by \ctn{Maitra15a} when only the number of individuals
tends to infinity, is clearly inapplicable.
Indeed, incorporation of random effects is asymptotically feasible only in our current asymptotic framework;
\ctn{Maitra15a} elucidate that inclusion of random effects does not make sense asymptotically unless the domains 
of observations are also increased indefinitely.
Also, only our current asymptotic framework allows different sets of time-dependent covariates for different individuals.

It is important to remark that the diffusion coefficient becomes known once the continuous process is completely
observed; see \ctn{Robert01}. Hence, following \ctn{Maitra15a} we assume that the diffusion coefficient is known,
and is not involved in the model selection problem.

We begin by establishing an asymptotic theory of Bayes factor for two
competing individual $SDE$'s, and then extend the theory to systems of $SDE$'s. In this context it 
is important to draw attention
to the fact that even this relatively simple problem of comparing any two individual $SDE$'s using Bayes factors
has not yet been considered in the literature.
Our investigation in this simpler case, however, faced with an apparently negative result; 
the associated Bayes factor failed to be consistent in the sense that the relevant variance failed
to converge to zero, even though convergence of the log-Bayes factor in expectation is ensured.
Despite this, we have been able to utilise this result to establish almost sure exponential convergence of the
Bayes factor when the number of individuals are also allowed to increase indefinitely.

The rest of our article is structured as follows. 
We begin with formalization of our set-up in Section \ref{sec:bf_sde}, while
we provide the necessary assumptions and results in Section \ref{sec:case2}. 
In Section \ref{sec:result_BF_consistency}
we investigate the asymptotics of Bayes factor for comparing two individual $SDE$'s.
We illustrate our results with a special case in Section \ref{sec:illustration}.
In Section \ref{sec:asymp_BF_nT} we exploit the asymptotic theory of Bayes factors developed for
comparing individual $SDE$'s to construct a convergence theory of Bayes factors comparing
systems of $SDE$'s in both $iid$ and non-$iid$ cases. 
In Section \ref{sec:simulated_data} we carry out two simulation studies to demonstrate
that Bayes factor yields the correct set of covariates in our $SDE$ models even in non-asymptotic
cases, and in Section \ref{sec:truedata}, we model a real, company-wise national stock exchange data set,
using a system of $SDE$'s, each consisting of a plausible set of covariates,
and obtain the best possible sets of covariate combinations
for the companies, using Bayes factor.
We summarize our contributions and provide concluding remarks in Section \ref{sec:conclusion}.

\section{Formalization of the model selection problem in the $SDE$ set-up when
$n\rightarrow\infty$ and $T_i\rightarrow\infty$ for every $i$}
\label{sec:bf_sde}

That the systems considered by us are well-defined and the exact likelihoods are computable, are guaranteed by
assumption (H2$^{\prime\prime}$) in Section \ref{sec:case2}. 
For our purpose we consider the filtration ($\mathcal F_t^W,t\geq 0$),
where $\mathcal F_t^W=\sigma(W_i(s),s\leq t)$. Each process $W_i$ is a $(\mathcal F_t^W , t\geq 0)$-adapted Brownian
motion.

Here we consider the set-up where, for $i=1,2,\ldots,n$,
\begin{equation}
d X_i(t)=\phi_{i,\bxi^{(i)}_0}(t)b_{\bbeta^{(i)}_{0}}(t,X_i(t))dt+\sigma(t,X_i(t))dW_i(t)
\label{eq:sde3}
\end{equation}
and
\begin{equation}
d X_i(t)=\phi_{i,\bxi^{(i)}_1}(t)b_{\bbeta^{(i)}_{1}}(t,X_i(t))dt+\sigma(t,X_i(t))dW_i(t),
\label{eq:sde4}
\end{equation}
where, $X_i(0)=x^i$ is the initial value of the stochastic process $X_i(t)$, 
which is assumed to be continuously observed on the time interval $[0,T_i]$; $T_i>0$.
We consider (\ref{eq:sde3}) as representing the true model and (\ref{eq:sde4}) is any other model.

It is useful to remark that we must analyze the same data set with respect to two different models for the purpose
of model selection. Hence, even though the distribution of the underlying stochastic process under the two models are different, 
for notational convenience we denote the process by $X_i(t)$ under both the models, relying on the context 
and the model-specific parameters to naturally clarify the distinction.


\subsection{Inclusion of time-dependent covariates}
\label{subsec:covariates}
We model $\phi_{i,\bxi^{(i)}_j}(t)$ for $j=0,1$, and $i=1,\ldots,n$, as  
\begin{equation}
\phi_{i,\bxi^{(i)}_j}(t)=\phi_{i,\bxi^{(i)}_j}(\bz_i(t))
=\xi^{(i)}_{0j}+\xi^{(i)}_{1j}g_1(z_{i1}(t))+\xi^{(i)}_{2j}g_2(z_{i2}(t))+\cdots+\xi^{(i)}_{pj}g_p(z_{ip}(t)),
\label{eq:phi_model}
\end{equation}
where 
$\bz_i(t)=(z_{i1}(t),z_{i2}(t),\ldots,z_{ip}(t))$ is the set of available covariate information 
corresponding to the $i$-th individual, depending upon time $t$. Following \ctn{Maitra15a} 
we assume $\bz_i(t)$ is continuous in $t$,
$z_{il}(t)\in \bZ_l$  where $\bZ_l$ is  compact and 
$g_l:\bZ_l\rightarrow \mathbb R$ is  continuous, for $l=1,\ldots,p$.
We let $\boldcal Z=\bZ_1\times\cdots\times\bZ_p$, and 
$\mathfrak Z=\left\{\bz(t)\in\boldcal Z:t\in[0,\infty)~\mbox{such that}~\bz(t)~\mbox{is continuous in}~t\right\}$.
Hence, $\bz_i\in\boldcal Z$ for all $i$.

\subsection{The random effects set-up}
\label{subsec:random_effects_2}

In (\ref{eq:sde3}), $\btheta^{(i)}_{0}=\left(\bbeta^{(i)}_{0},\xi^{(i)}_{00},\xi^{(i)}_{10},\ldots,\xi^{(i)}_{p0}
\right)=\left(\bbeta^{(i)}_{0},\bxi_0^{(i)}\right)$ stands for 
the true parameters, and  
$\btheta^{(i)}_{1}=\left(\bbeta^{(i)}_{1},\xi^{(i)}_{01},\xi^{(i)}_{11},\ldots,\xi^{(i)}_{p1}\right)
=\left(\bbeta^{(i)}_{1},\bxi_1^{(i)}\right)$ are
the parameters associated with (\ref{eq:sde4}). 
Let $\btheta^{(i)}_j\in\bTheta=\mathfrak B\times\bGamma$ for all $i$,
where both $\mathfrak B$ and $\bGamma$ are compact spaces. 
We also assume that for $i=1,2,\ldots,n$,
$$\btheta^{(i)}_{1}\stackrel{iid}{\sim}\pi,$$
where $\pi$ is some specified distribution on $\bTheta$.

Hence, the above describes a random effects set-up.
Observe that if $\xi^{(i)}_{lj}=0$ for $l=1,\ldots,p$, and for $i=1,\ldots,n$, then it reduces
to the random effects model of \ctn{Maud12}, showing that the latter is a special case of our model.

As is well-known, even though the term ``prior" is not appropiate
for the random effects coefficients, operationally there is no difference between
a prior and a distribution for random effects in the Bayesian paradigm. Somewhat abusing the terminology, we continue
to refer to the distribution of the $iid$ random effects coeffcients, $\pi$, as the relevant prior.

\subsection{Covariate and drift function selection}
\label{subsec:model_selection_drift_covariates}

The key difference between our current model selection idea and that of \ctn{Maitra15a} is that
here, for every individual, there is an independent model selection problem. In other words, for each $i$,
one needs to choose between $\btheta^{(i)}_0$ and $\btheta^{(i)}_1$. This involves selection of
perhaps different sets of covariates for different $i$ with respect to the coefficients
$\bxi^{(i)}_j$, and different drift functions $b_{\bbeta^{(i)}_j}$. 
Obviously, the dimensions of $\bxi^{(i)}_0$ and $\bxi^{(i)}_1$ are allowed to differ for each $i$;
likewise, for every $i$, the dimensions of $\bbeta^{(i)}_0$ and $\bbeta^{(i)}_1$ may be different as well.
Thus, from this perspective, our current model selection framework appears to be more general compared to that 
of \ctn{Maitra15a}, who consider
the same set of parameters $\bxi_j$ and $\bbeta_j$ for all the individuals, allowing only a fixed set of covariates
for every subject.

%

\subsection{Form of the Bayes factor in our set-up}
\label{subsec:bf_form}
For $j=0,1$, we
first define the following quantities:
\begin{equation}
U_{i,\btheta^{(i)}_j} =\int_0^{T_i}\frac{\phi_{i,\bxi^{(i)}_j}(s)b_{\bbeta^{(i)}_j}(s,X_i(s))}{\sigma^2(s,X_i(s))}dX_i(s),
\quad\quad V_{i,\btheta^{(i)}_j} =\int_0^{T_i}\frac{\phi^2_{i,\bxi^{(i)}_j}(s)b_{\bbeta^{(i)}_j}^2(s,X_i(s))}
{\sigma^2(s,X_i(s))}ds
\label{eq:u_v}
\end{equation}
for $j=0,1$ and $i=1,\ldots,n$. 

Let $\bC_{T_i}$ denote the space of real continuous functions $(x(t), t \in [0,T_i ])$ defined on $[0,T_i ]$, 
endowed with the $\sigma$-field $\mathcal C_{T_i}$ associated with the topology of uniform convergence 
on $[0,T_i ]$. We consider the distribution $P^{x_i,T_i,\bz_i}_j$ on $(C_{T_i} ,\mathcal C_{T_i})$ 
of $(X_i (t), t\in [0,T_i])$ given by (\ref{eq:sde3}) and (\ref{eq:sde4}) for $j=0,1$. 
We choose the dominating measure $P_i $ as the distribution of (\ref{eq:sde3}) and (\ref{eq:sde4}) with null
drift. So, for $j=0,1$, 
\begin{align}
\frac{dP^{x_i,T_i,\bz_i}_j}{dP_i}=f_{i,\btheta^{(i)}_j}(X_i)
&=\exp\left(U_{i,\btheta^{(i)}_j}-\frac{V_{i,\btheta^{(i)}_j}}{2}\right), 
\label{eq:densities}
\end{align}
where $f_{i,\btheta^{(i)}_0}(X_i)$ denotes the true density and $f_{i,\btheta^{(i)}_1}(X_i)$ stands for the other density
associated with the modeled $SDE$.

For each $i=1,\ldots,n$, letting $X_{i,a,b}$ denote the $i$-th process 
observed on $[a,b]$ for any $0\leq a<b<\infty$,
\begin{equation}
I_{x^i,T_i,\bz_i}=\int_{\bTheta} \frac{f_{i,\btheta^{(i)}_1}(X_{i,0,T_i})}
{f_{i,\btheta^{(i)}_0}(X_{i,0,T_i})}\pi\left(\btheta^{(i)}_{1}\right)d\btheta^{(i)}_{1}
\label{eq:I_iT}
\end{equation}
denotes the Bayes factor associated with the $i$-th equation of the above two systems of equations. 
Assuming that the $SDE$'s (\ref{eq:sde3}) and (\ref{eq:sde4}) are independent for $i=1,\ldots,n$, 
\begin{equation*}
I_{n,T_1,\ldots,T_n}=\prod_{i=1}^nI_{x^i,T_i,\bz_i}
\label{eq:I_nT}
\end{equation*}
is the Bayes factor comparing the entire systems of $SDE$'s (\ref{eq:sde3}) and (\ref{eq:sde4}). 

Comparisons between a collection of different models using Bayes factor, none of which may be the true model, is 
expected to favour that model which minimizes the Kullback-Leibler
divergence from the true model.

\subsection{The $iid$ and the non-$iid$ cases}
\label{subsec:iid_non_iid_2}
We are interested in studying the properties of $I_{n,T_1,\ldots,T_n}$ in both $iid$ and non-$iid$ cases
when $n\rightarrow\infty$ and $T_i\rightarrow\infty$. 
In the $iid$ set-up, we assume that $x^i=x$, $T_i=T$ and $\btheta^{(i)}_j=\left(\bbeta_j^{(i)},\xi^{(i)}_{0j}\right)$, 
for $i=1,\ldots,n$ and $j=0,1$. In the non-$iid$ case we relax these assumptions.
However, for simplicity, we assume
$T_i=T$ for each $i$, even in the non-$iid$ set-up, so that in our asymptotic framework we study convergence of 
\begin{equation}
\tilde I_{n,T}=\prod_{i=1}^nI_{i,T},
\label{eq:tilde_I_nT}
\end{equation}
as $n\rightarrow\infty$ and $T\rightarrow\infty$, where $I_{i,T}=I_{x^i,T,\bz_i}$.

\subsection{A key relation between $U_{i,\btheta^{(i)}_j}$ and $V_{i,\btheta^{(i)}_j}$ in the context
of model selection using Bayes factors}
\label{subsec:key_relation_U_V}
An useful relation between $U_{i,\btheta^{(i)}_j}$ and $V_{i,\btheta^{(i)}_j}$ which we will often make use of
in this paper is as follows.
\begin{align}
U_{i,\btheta^{(i)}_j}&=\int_0^{T_i}\frac{\phi_{i,\bxi^{(i)}_j}(s)b_{\bbeta^{(i)}_j}\left(X_i(s)\right)}
{\sigma^2\left(X_i(s)\right)}dX_i(s)\notag\\
&=\int_0^{T_i}\frac{\phi_{i,\bxi^{(i)}_j}(s)b_{\bbeta^{(i)}_j}\left(X_i(s)\right)}{\sigma^2\left(X_i(s)\right)}
\left[\phi_{i,\bxi^{(i)}_0}(s)b_{\bbeta^{(i)}_0}\left(X_i(s)\right)ds+\sigma\left(X_i(s)\right)dW_i(s)\right]\notag\\
&=\int_0^{T_i}\frac{\phi_{i,\bxi^{(i)}_j}(s)\phi_{i,\bxi^{(i)}_0}(s)b_{\bbeta^{(i)}_j}
\left(X_i(s)\right)b_{\bbeta^{(i)}_0}\left(X_i(s)\right)}
{\sigma^2\left(X_i(s)\right)}ds
+\int_0^{T_i}\frac{\phi_{i,\bxi^{(i)}_j}(s)b_{\bbeta^{(i)}_j}\left(X_i(s)\right)}
{\sigma\left(X_i(s)\right)}dW_i(s)\notag\\
&=V_{i,\btheta^{(i)}_0,\btheta^{(i)}_j}
+\int_0^{T_i}\frac{\phi_{i,\bxi^{(i)}_j}(s)b_{\bbeta^{(i)}_j}\left(X_i(s)\right)}
{\sigma\left(X_i(s)\right)}dW_i(s),
\label{eq:u_v_relation}
\end{align}
with
\begin{equation}
V_{i,\btheta^{(i)}_0,\btheta^{(i)}_j}=
\int_0^{T_i}\frac{\phi_{i,\bxi^{(i)}_j}(s)\phi_{i,\bxi^{(i)}_0}(s)b_{\bbeta^{(i)}_j}
\left(X_i(s)\right)b_{\bbeta^{(i)}_0}
\left(X_i(s)\right)}{\sigma^2\left(X_i(s)\right)}ds.
\label{eq:V_0_j}
\end{equation}
Note that $V_{i,\btheta^{(i)}_0}=V_{i,\btheta^{(i)}_0,\btheta^{(i)}_0}$ and 
$V_{i,\btheta^{(i)}_1}=V_{i,\btheta^{(i)}_1,\btheta^{(i)}_1}$.
Also note that, for $j=0,1$, for each $i$,
\begin{equation}
E_{\btheta^{(i)}_0}\left[\int_0^{T_i}\frac{\phi_{i,\bxi^{(i)}_j}(s)b_{\bbeta^{(i)}_j}\left(X_i(s)\right)}
{\sigma\left(X_i(s)\right)}dW_i(s)\right]=0,
\label{eq:zero_mean}
\end{equation}
so that $E_{\btheta^{(i)}_0}\left(U_{i,\btheta^{(i)}_j}\right)
=E_{\btheta^{(i)}_0}\left(V_{i,\btheta^{(i)}_0,\btheta^{(i)}_j}\right)$.


\section{Requisite assumptions and results for the asymptotic theory of Bayes factor 
when $n\rightarrow\infty$ and $T\rightarrow\infty$}
\label{sec:case2}

%

All our following assumptions and results are true for each $i$, in particular true for each 
$\bbeta_j^{(i)},\bxi_j^{(i)}$ and consequently for $U_{i,\btheta_j^{(i)}}, 
V_{i,\btheta_j^{(i)}},V_{i,\btheta_0^{(i)},\btheta_j^{(i)}}$. For the sake of notational simplicity 
we provide all the assumptions and results without mentioning $i$ at every stage. 
We make the following assumptions:
\begin{itemize}
\item[(H1$^{\prime\prime}$)] The parameter space $\bTheta=\mathfrak B\times\bGamma$ such that 
$\bGamma$ and $\mathfrak B$ are compact.
\end{itemize}

\begin{itemize}
\item[(H2$^{\prime\prime}$)] 
For $j=0,1$, given any $s$, $\bbeta_j$,
$b_{\bbeta_j}(s,\cdot)$,  $\sigma(s,\cdot)$ are $C^1$ on $\mathbb R$;  
we also assume that $b^2_{\bbeta_j}(s,x)\leq K_1(1+x^2+\|\bbeta_j\|^2)$  
and $\sigma^2(x)\leq K_2(1+x^2)$
for all $s\in [0,T]$, $x\in\mathbb R$, for some $K_1,K_2>0$. 
By (H1$^{\prime\prime}$) it follows as before that for $s\in[0,T]$, 
$b^2_{\bbeta_j}(s,x)\leq K(1+x^2)$  and $\sigma^2(s,x)\leq K(1+x^2)$
for all $x\in\mathbb R$, for some $K>0$. 

\end{itemize}
Because of (H2$^{\prime\prime}$) it follows from Theorem 4.4 of \ctn{Mao11}, page 61, that  
for all $T>0$, and any $k\geq 2$, 
\begin{equation}
E\left(\underset{s\in [0,T]}{\sup}~|X_i(s)|^k\right)\leq\left(1+3^{k-1}E|X_i(0)|^k\right)\exp\left(\tilde\vartheta T\right),
\label{eq:moment1}
\end{equation}
where
$$\tilde\vartheta=\frac{1}{6}\left(18K\right)^{\frac{k}{2}}T^{\frac{k-2}{2}}\left[T^{\frac{k}{2}}
+\left(\frac{k^3}{2(k-1)}\right)^{\frac{k}{2}}\right].$$

Specifically, for any $k\geq 2$, we can write, as $T\rightarrow\infty$,
\begin{equation}
E\left(\underset{s\in [0,T]}{\sup}~|X(s)|^k\right)=o\left(\exp\left\{T^{k+1}\right\}\right).
\label{eq:moment1_T}
\end{equation}
%
%
We further assume the following conditions.
\begin{itemize}
\item[(H3$^{\prime\prime}$)] $b_{\bbeta_j}(s,x)$ is continuous 
in $(x,\bbeta_j)$.

\item[(H4$^{\prime\prime}$)] 
For $s\in[0,T]$ and $j=0,1$, 
$\frac{b_{\bbeta_j}^2(s,x)}{\sigma^2(s,x)}$ and $\frac{b_{\bbeta_1}(s,x)b_{\bbeta_0}(s,x)}{\sigma^2(s,x)}$ 
satisfy the following:
\begin{align}
&\kappa_j(\bbeta_j) + \frac{K_{\bbeta_j}\left(1+x^2+\|\bbeta_j\|^2\right)}{c_j+\exp\left(T^5\right)}
<\frac{b_{\bbeta_j}^2(s,x)}{\sigma^2(s,x)}\notag\\
&\quad\quad\quad<\kappa_j(\bbeta_j)+
\frac{M_{\bbeta_j}\left(1+x^2+\|\bbeta_j\|^2\right)}{d_j+\exp\left(T^5\right)},
\label{eq:H4_prime_prime_1}
\end{align}
and
\begin{align}
&\bar{\kappa}(\bbeta_0,\bbeta_1) + \frac{K_{\bbeta_0,\bbeta_1}
\left(1+x^2+\|\bbeta_0\|^2+\|\bbeta_1\|^2\right)}{\bar c+\exp\left(T^5\right)}< 
\frac{b_{\bbeta_1}(s,x)b_{\bbeta_0}(s,x)}{\sigma^2(s,x)}\notag\\
&\quad\quad\quad<\bar{\kappa}(\bbeta_0,\bbeta_1)+\frac{M_{\bbeta_0,\bbeta_1}
\left(1+x^2+\|\bbeta_0\|^2+\|\bbeta_1\|^2\right)}{\bar d+\exp\left(T^5\right)},
\label{eq:H4_prime_prime_2}
\end{align}
where $0<c_j,d_j,\bar c,\bar d<\infty$, are some constants;
$\kappa_j(\bbeta_j)$ are positive, continuous functions of $\bbeta_j$; $\bar{\kappa}(\bbeta_0,\bbeta_1)$
is a continuous function of $(\bbeta_0,\bbeta_1)$;
$K_{\bbeta_j}$, $M_{\bbeta_j}$ are continuous in $\bbeta_j$, for $j=0,1$,
and $K_{\bbeta_0,\bbeta_1}$, $M_{\bbeta_0,\bbeta_1}$ are continuous in $(\bbeta_0,\bbeta_1)$.

\item[(H5$^{\prime\prime}$)]

(i) We assume that $\boldcal Z= \bZ_1\times\bZ_2\times\cdots\times\bZ_p$ is the space of the covariates
where $\bZ_l$ is compact for $l=1,\ldots,p$, and for every $t\geq 0$, 
$\bz_i(t)=(z_{i1}(t),z_{i2}(t),\ldots,z_{ip}(t))\in\boldcal Z$ for $i=1,\ldots,n$.
Also, we assume that $\bz_i(t)$ are continuous in $t$ for every $i$, so that $\bz_i\in\mathfrak Z$, for every $i$.

(ii) For $j=0,1$, and for $t\geq 0$, we assume that the vector of covariates $\bz_i(t)$ is related to the $i$-th $SDE$
of the $j$-th model via
$$\phi_{i,\bxi^{(i)}_j}(t)=\phi_{i,\bxi^{(i)}_j}(\bz_i(t))=\xi^{(i)}_{0j}+\sum_{l=1}^p\xi^{(i)}_{lj}g_l(\bz_i(t)),$$
where, for $l=1,\ldots,p$, $g_l:\bZ_l\rightarrow \mathbb R$ is 
continuous. Notationally, when reference to the $i$-th individual is self-explanatory, 
we shall denote the function $\xi_{0j}+\sum_{l=1}^p\xi_{lj}g_l$
by $\phi_{\bxi_j}$.

(iii) For $l=1,2,\ldots,p$, for $i=1,\ldots,n$, and for $t\geq 0$,
\begin{equation}
\frac{1}{n}\sum_{i=1}^n g_l(z_{il}(t))\rightarrow c_{l}(t),
\label{eq:H5_prime_1_T}
\end{equation}
and
\begin{equation}
\frac{1}{n}\sum_{i=1}^n g_l(z_{il}(t))g_m(z_{im}(t))\rightarrow c_l(t)c_m(t),
\label{eq:H5_prime_2_T}
\end{equation}
as $n\rightarrow\infty$, where $\left\{c_l(t):t\geq 0\right\}$ are real constants for $l=1,\ldots,p$.

(iv) For $l=1,2,\ldots,p$, and for $i=1,\ldots,n$,
\begin{equation}
\underset{T\rightarrow\infty}{\lim}~\frac{1}{T}\int_0^T g_l(z_{il}(s))ds =\bar c^{(1)}_{il}\label{eq::H5_prime_3_T}
\end{equation}
and
\begin{equation}
\underset{T\rightarrow\infty}{\lim}~\frac{1}{T}\int_0^T g_l(z_{il}(s))g_m(z_{im}(s))ds
=\bar c^{(2)}_{ilm},\label{eq::H5_prime_4_T}
\end{equation}
where $\bar c^{(1)}_{il}$ and $\bar c^{(2)}_{ilm}$ are real constants.
\end{itemize}
\begin{remark}
\label{remark:time_dependent_covariates2}
Observe that although (H4$^{\prime\prime}$) is seemingly restrictive in the sense that the ratios
$\frac{b_{\bbeta_j}^2(s,x)}{\sigma^2(s,x)}$ and $\frac{b_{\bbeta_1}(s,x)b_{\bbeta_0}(s,x)}{\sigma^2(s,x)}$
are approximately independent of the underlying stochastic process, assumption (H5$^{\prime\prime}$)
attempts to compensate for the restrictions by providing a rich structure to $\phi_{\bxi_j}$ consisting of 
covariate information varying continuously with time. 
Hence, assumption (H4$^{\prime\prime}$) need not be viewed as restrictive. 
\end{remark}

\ctn{Maitra15a} argue that (\ref{eq:H5_prime_1_T}) and (\ref{eq:H5_prime_2_T}) hold
if one assumes that for $i=1,\ldots,n$, and  
$l=1,\ldots,p$, the covariates $z_{il}$ are observed realizations
of stochastic processes that are $iid$ for $i=1,\ldots,n$, for all $l=1,\ldots,p$, and that 
for $l\neq m$, the processes generating $z_{il}$ and $z_{im}$ are independent. 
In other words, although we assume the covariates to be non-random, in
essence, it may be assumed $g_l(z_{il}(t))$ and $g_m(z_{im}(t))$ are  
uncorrelated for $l\neq m$. 

In order that (H5$^{\prime\prime}$) (iv) holds, one needs to further assume that
the relevant stochastic processes converge to appropriate stationary distributions.
For example, $z_{il}(t)$ may be realizations of Markov processes which are
irreducible (with respect to some appropriate measure), aperiodic, positive recurrent and possses invariant
distributions; see, for example, \ctn{Konto03}. 

It follows from (H5$^{\prime\prime}$) (iv), that, 
\begin{align}
&\underset{T\rightarrow\infty}{\lim}~ \frac{1}{T}\int_0^T\phi_{i,\bxi^{(i)}_j}(\bz_i(s))ds\notag\\
&= \xi^{(i)}_{0j}+\sum_{l=1}^p\xi^{(i)}_{lj}\bar c^{(1)}_{il}\notag\\
&=\bar\phi^{(1)}_{i,\bxi^{(i)}_j}~\mbox{(say)},
\label{eq:phi_conv}
\end{align}
\begin{align}
&\underset{T\rightarrow\infty}{\lim}~\frac{1}{T}\int_0^T\phi^2_{i,\bxi^{(i)}_j}(\bz_i(s))ds\notag\\
&= \left\{\xi^{(i)}_{0j}\right\}^2+2\xi^{(i)}_{0j}\sum_{l=1}^p\xi^{(i)}_{lj}\bar c^{(1)}_{il}
+\sum_{l=1}^p\sum_{m=1}^p\xi^{(i)}_{lj}\xi^{(i)}_{mj}\bar c^{(2)}_{ilm}\notag\\
&=\bar\phi^{(2)}_{i,\bxi^{(i)}_j},~\mbox{(say)},
\label{eq:phisq_conv}
\end{align}
and
\begin{align}
&\underset{T\rightarrow\infty}{\lim}~\frac{1}{T}\int_0^T\phi_{i,\bxi^{(i)}_0}(\bz_i(s))\phi_{i,\bxi^{(i)}_1}(\bz_i(s))ds\notag\\
&=\xi^{(i)}_{00}\xi^{(i)}_{01}+\xi^{(i)}_{00}\sum_{l=1}^p\xi^{(i)}_{l1}\bar c^{(1)}_{il}
+\xi^{(i)}_{01}\sum_{l=1}^p\xi^{(i)}_{l0}\bar c^{(1)}_{il}+
\sum_{l=1}^p\sum_{m=1}^p\xi^{(i)}_{l0}\xi^{(i)}_{m1}\bar c^{(2)}_{ilm}\notag\\
&=\bar\phi^{(2)}_{i,\bxi^{(i)}_0,\bxi^{(i)}_1},~\mbox{(say)}.
\label{eq:phiphi_conv}
\end{align}
When $i$ is clear from the context, we shall often use the notations 
$\bar\phi^{(1)}_{\bxi_j}$, $\bar\phi^{(2)}_{\bxi_j}$ and $\bar\phi^{(2)}_{\bxi_0,\bxi_1}$.

Note that, (\ref{eq:phi_conv}), (\ref{eq:phisq_conv}) and (\ref{eq:phiphi_conv}) are limits of expectations
with respect to the uniform distribution on $[0,T]$. Hence, by the Cauchy-Schwartz inequality it follows that
\begin{equation}
\bar\phi^{(2)}_{\bxi_0,\bxi_1}
\leq \sqrt{\bar\phi^{(2)}_{\bxi_0}\bar\phi^{(2)}_{\bxi_1}}.
\label{eq:phi_cs}
\end{equation}

The following lemmas will be useful in our proceedings. The proofs of these lemmas are provided in sections S-1, S-2 and S-3 of the supplement. 
\begin{lemma}
\label{lemma:uniform_convergence}
The limits
$\bar\phi^{(1)}_{\bxi_1}$, $\bar\phi^{(2)}_{\bxi_1}$ and $\bar\phi^{(2)}_{\bxi_0,\bxi_1}$
are continuous in $\bxi_1$.
\end{lemma}

\begin{lemma}
\label{lemma:convergence1}
Assume (H1$^{\prime\prime}$) -- (H5$^{\prime\prime}$). Then, the following hold:
\begin{align}
(i)\quad &E_{\btheta_0}\left(\frac{V_{\btheta_j,T}}{T}\right)\rightarrow 
\bar \phi^{(2)}_{\bxi_j}\kappa_j(\bbeta_j);j=0,1.\label{eq:E_V_convergence}\\
(ii)\quad &E_{\btheta_0}\left(\frac{V_{\btheta_0,\btheta_1,T}}{T}\right)
\rightarrow  \bar \phi^{(2)}_{\bxi_0,\bxi_1}\bar{\kappa}(\bbeta_0,\bbeta_1).
\label{eq:E_V_convergence2}\\
(iii)\quad &E_{\btheta_0}\left(\frac{U_{\btheta_0,T}}{T}\right)\rightarrow 
\bar\phi^{(2)}_{\bxi_0}\kappa_0(\bbeta_0).\label{eq:E_U_convergence_0}\\
(iv)\quad &E_{\btheta_0}\left(\frac{U_{\btheta_1,T}}{T}\right)\rightarrow 
\bar\phi^{(2)}_{\bxi_0,\bxi_1}\bar{\kappa}(\bbeta_0,\bbeta_1).\label{eq:E_U_convergence}\\
(v)\quad &\frac{1}{T}\int_0^T\frac{\phi_{\bxi_j}(s)b_{\bbeta_j}(s,X(s))}{\sigma(s,X(s))}dW(s)
\stackrel{a.s.}{\longrightarrow} 0; j=0,1,
\label{eq:integral_convergence}\\
(vi)\quad &\frac{V_{\btheta_j,T}}{T}\stackrel{a.s.}{\longrightarrow} 
 \bar \phi^{(2)}_{\bxi_j}\kappa_j(\bbeta_j); j=0,1, \label{eq:V_convergence}\\
(vii)\quad &\frac{V_{\btheta_0,\btheta_1,T}}{T}\stackrel{a.s.}{\longrightarrow} 
\bar \phi^{(2)}_{\bxi_0,\bxi_1}\bar{\kappa}(\bbeta_0,\bbeta_1), 
\label{eq:V_convergence2}\\
(viii)\quad &\frac{U_{\btheta_0,T}}{T}\stackrel{a.s.}{\longrightarrow} 
\bar\phi^{(2)}_{\bxi_0}\kappa_0(\bbeta_0),
\label{eq:U_convergence_0}\\
(ix)\quad &\frac{U_{\btheta_1,T}}{T}\stackrel{a.s.}{\longrightarrow} 
\bar\phi^{(2)}_{\bxi_0,\bxi_1}\bar{\kappa}(\bbeta_0,\bbeta_1),
\label{eq:U_convergence}
\end{align}
In the above, $``\stackrel{a.s.}{\longrightarrow}"$ denotes convergence ``almost surely" as $T\rightarrow\infty$
with respect to $X$ (under $\btheta_0$), and the expectations are also with respect to $X$ (under $\btheta_0$).
\end{lemma}

\begin{lemma}
\label{lemma:phi_kappa_relations}
Assume (H1$^{\prime\prime}$) -- (H5$^{\prime\prime}$). Then, the following holds:
\begin{equation}
\bar\phi^{(2)}_{\bxi_0,\bxi_1}\bar\kappa(\bbeta_0,\bbeta_1)
\leq \sqrt{\bar\phi^{(2)}_{\bxi_0}\bar\phi^{(2)}_{\bxi_1}}
\times\sqrt{\kappa_0(\bbeta_0)\kappa_1(\bbeta_1)}.
\label{eq:phi_kappa_relations}
\end{equation}
\end{lemma}

\section{Convergence of Bayes factor with respect to time when two individual $SDE$'s are compared}
\label{sec:result_BF_consistency}

From the system of $SDE$'s defined by (\ref{eq:sde3}) and (\ref{eq:sde4}) we now consider the $i$-th individual only. 
To avoid notational complexity we denote $X_i$ simply by $X$. Consequently, $\phi_{i,\bxi_j}(t)$ and $T_i$ 
will be denoted by $\phi_{\bxi_j}(t)$ and $T$, respectively. In connection with the $i$-th individual 
we consider the following two $SDE$'s:
\begin{equation}
d X(t)=\phi_{\bxi_0}(t)b_{\bbeta_0}(t,X(t))dt+\sigma(t,X(t))dW(t)
\label{eq:sde1_T}
\end{equation}
and
\begin{equation}
d X(t)=\phi_{\bxi_1}(t)b_{\bbeta_1}(t,X(t))dt+\sigma(t,X(t))dW(t).
\label{eq:sde2_T}
\end{equation}

For any $t\in[0,T]$, for $j=0,1$, let
\begin{align}
U_{\btheta_j,t} &=\int_0^t\frac{\phi_{\bxi_j}(s)b_{\bbeta_j}(s,X(s))}{\sigma^2(s,X(s))}dX(s),\quad 
V_{\btheta_j,t} =\int_0^t\frac{\phi^2_{\bxi_j}(s)b^2_{\bbeta_j}(s,X(s))}{\sigma^2(s,X(s))}ds,\notag\\
V_{\btheta_0,\btheta_j,t} 
&=\int_0^t\frac{\phi_{\bxi_j}(s)b_{\bbeta_j}(s,X(s))\phi_{\bxi_0}(s)b_{\bbeta_0}(s,X(s))}{\sigma^2(s,X(s))}ds.
\label{eq:u_v_t}
\end{align}
Note that $V_{\btheta_0,t}= V_{\btheta_0,\btheta_0,t}$ and 
$V_{\btheta_1,t}= V_{\btheta_1,\btheta_1,t}$.
We also let
\begin{equation}
f_{\btheta_j,t}(X_{0,t})=\exp\left(U_{\btheta_j,t}-\frac{V_{\btheta_j,t}}{2}\right).
\label{eq:f_t}
\end{equation}
Here we are interested in asymptotic properties of the Bayes factor, given by
\begin{equation}
I_T=\int \frac {f_{\btheta_1,T}(X_{0,T})}{f_{\btheta_0,T}(X_{0,T})}\pi (d\btheta_1),
\label{eq:bf_T}
\end{equation}
as $T\rightarrow\infty$.

For our purpose, let us define, for any $h>0$,
\begin{align}
U_{\btheta_j,t,t+h} &=\int_t^{t+h}\frac{\phi_{\bxi_j}(s)b_{\bbeta_j}(s,X(s))}{\sigma^2(s,X(s))}dX(s),\quad 
V_{\btheta_j,t,t+h} =\int_t^{t+h}\frac{\phi^2_{\bxi_j}(s)b^2_{\bbeta_j}(s,X(s))}{\sigma^2(s,X(s))}ds,\notag\\
V_{\btheta_0,\btheta_j,t,t+h} 
&=\int_t^{t+h}\frac{\phi_{\bxi_j}(s)b_{\bbeta_j}(s,X(s))\phi_{\bxi_0}(s)b_{\bbeta_0}(s,X(s))}{\sigma^2(s,X(s))}ds.
\label{eq:u_v_t2}
\end{align}
Observe, as before, that $V_{\btheta_0,t,t+h}= V_{\btheta_0,\btheta_0,t,t+h}$ 
and $V_{\btheta_1,t,t+h}= V_{\btheta_1,\btheta_1,t,t+h}$.
We let
\begin{equation}
f_{\btheta_j,t,t+h}(X_{t,t+h})=\exp\left(U_{\btheta_j,t,t+h}-\frac{V_{\btheta_j,t,t+h}}{2}\right),
\label{eq:f_t2}
\end{equation}
where, for any $0\leq a<b<\infty$, $X_{a,b}$ denotes a path of the process $X$ from $a$ to $b$.
For any $t>0$ and $h>0$, we define 
\begin{equation}
\tilde{\mathcal K}(f_{\btheta_0,t,t+h},f_{\btheta_1,t,t+h})
=E_{\btheta_0}\left[\log\frac{f_{\btheta_0,t,t+h}}{f_{\btheta_1,t,t+h}}\right],
\label{eq:kl1_T}
\end{equation}
where $E_{\btheta_0}\equiv E_{f_{\btheta_0,t+h}}$.
Note that although the expectation is with respect to $f_{\btheta_0,t+h}$, which is not the same as $f_{\btheta_0,t,t+h}$,
(\ref{eq:kl1_T}) is still the Kullback-Leibler divergence between $f_{\btheta_0,t,t+h}$ and $f_{\btheta_1,t,t+h}$.
Also since in our case, for $j=0,1$,
\begin{align}
f_{\btheta_j,t,t+h}\left(X_{t,t+h}\right)&=\exp\left(U_{\btheta_j,t,t+h}-\frac{V_{\btheta_j,t,t+h}}{2}\right)\notag\\
&=\exp\left((U_{\btheta_j,t+h}-U_{\btheta_j,t})
-\frac{(V_{\btheta_j,t+h}-V_{\btheta_j,t})}{2}\right),
\label{eq:f_t3}
\end{align}
it follows that
\begin{equation}
\tilde{\mathcal K}(f_{\btheta_0,t,t+h},f_{\btheta_1,t,t+h})
=\mathcal K(f_{\btheta_0,t+h},f_{\btheta_1,t+h})-\mathcal K(f_{\btheta_0,t},f_{\btheta_1,t}),
\label{eq:kl2}
\end{equation}
where $\mathcal K(f_{\btheta_0,t+h},f_{\btheta_1,t+h})$ and $\mathcal K(f_{\btheta_0,t},f_{\btheta_1,t})$
are proper Kullback-Leibler divergences between $f_{\btheta_0,t+h}$, $f_{\btheta_1,t+h}$,
and $f_{\btheta_0,t}$, $f_{\btheta_1,t}$, respectively.
We now define 
\begin{align}
\tilde{\mathcal K}'_t(f_{\btheta_0},f_{\btheta_1})&=\lim_{h\rightarrow 0} 
\frac{\tilde {\mathcal K}(f_{\btheta_0,t,t+h},f_{\btheta_1,t,t+h})}{h}\notag\\
&=\frac{1}{2}\frac{d}{dt}E_{\theta_0}(V_{\btheta_0,t})
-\frac{d}{dt}E_{\theta_0}(V_{\btheta_0,\btheta_1,t})
+\frac{1}{2}\frac{d}{dt}E_{\theta_0}(V_{\btheta_1,t}).
\label{eq:kl3}
\end{align}
The expression (\ref{eq:kl3}) easily follows using (\ref{eq:f_t3}), 
the relation (\ref{eq:u_v_relation}) and (\ref{eq:zero_mean}).

\subsection{Pseudo  Kullback-Leibler $(\delta)$ property}
\label{subsec:pseudo_kl}
We make the following assumption: 
\begin{itemize}
\item[(H6$^{\prime\prime}$)]
For a fixed $\delta\geq 0$, the prior $\pi$ satisfies
\begin{equation}
\pi\left(\btheta_1\in\bTheta:\underset{t}{\inf}~\tilde{\mathcal K}'_t(f_{\btheta_0},f_{\btheta_1})\geq\delta\right)=1.
\label{eq:pi_prob_1}
\end{equation}
\end{itemize}

Let us define
\begin{equation}
\bar{\mathcal K}^{\infty}(f_{\btheta_0},f_{\btheta_1})
=\underset{T\rightarrow\infty}{\lim}~\frac{1}{T}\int_0^T\tilde{\mathcal K}'_t(f_{\btheta_0},f_{\btheta_1})dt.
\label{eq:kl_average_infinite}
\end{equation}
We assume the following:
\begin{itemize}
\item[(H7$^{\prime\prime}$)]
Given $\delta$ associated with (H6$^{\prime\prime}$), for any $c\geq 0$, the prior $\pi$ satisfies
\begin{equation}
\pi\left(\btheta_1\in\bTheta:\delta\leq\bar{\mathcal K}^{\infty}(f_{\btheta_0},f_{\btheta_1})\leq\delta+c\right)>0.
\label{eq:kl_property_T}
\end{equation}
\end{itemize}
We refer to property (H7$^{\prime\prime}$) as the pseudo Kullback-Leibler ($\delta$) property of the prior $\pi$.
Note that, (\ref{eq:kl3}), (\ref{eq:E_V_convergence}) and (\ref{eq:E_V_convergence2}) imply
\begin{align}
\bar{\mathcal K}^{\infty}(f_{\btheta_0},f_{\btheta_1})
&=\frac{\bar\phi^{(2)}_{\bxi_0}}{2}\kappa_0(\bbeta_0)
-\bar\phi^{(2)}_{\bxi_0,\bxi_1}\bar{\kappa}(\bbeta_0,\bbeta_1)
+\frac{\bar\phi^{(2)}_{\bxi_1}}{2}\kappa_1(\bbeta_1)\label{eq:kl_average_infinite2}\\
&\geq\frac{1}{2}\left(\sqrt{\bar\phi^{(2)}_{\bxi_0}\kappa_0(\bbeta_0)}
-\sqrt{\bar\phi^{(2)}_{\bxi_1}\kappa_1(\bbeta_1)}\right)^2\notag\\
&\geq 0,\notag
\end{align}
by Lemma \ref{lemma:phi_kappa_relations}.
Provided that (\ref{eq:pi_prob_1}) holds and the prior $\pi$ is dominated by the Lebesgue measure, 
the pseudo Kullback-Leibler ($\delta$)
property holds because of continuity of (\ref{eq:kl_average_infinite2}) 
in $\btheta_1=(\bbeta_1,\bxi_1)$ ensured by Lemma \ref{lemma:uniform_convergence}.

\subsection{$Q^*$ property}
\label{subsec:Q_star_T}
%
%
For $t\geq 0$, let $\mathcal F_{t}$ be the $\sigma$-algebra generated by $X(0)$ and the history of the process upto 
(and including) time $t$, and let $\pi_t(\btheta_1) =\pi(\btheta_1|\mathcal F_t)$ be the posterior
of $\btheta_1$ given $\mathcal F_t$. Also, let
\begin{align}
&\hat f_{t-h,t}\left(X_{t-h,t}\right)\notag\\
&=\int_{\btheta_1\in\bTheta}
\exp\left(\int_{t-h}^t\frac{\phi_{\bxi_1}(s)b_{\bbeta_1}(s,X(s))}{\sigma^2(s,X(s))}dX(s)
-\frac{1}{2}\int_{t-h}^t\frac{\phi_{\bxi_1}^2(s)b_{\bbeta_1}^2(s,X(s))}{\sigma^2(s,X(s))}ds\right)
\pi_{t-h}(d\btheta_1)\notag\\
&=E\left[\exp\left(\int_{t-h}^t\frac{\phi_{\bxi_1}(s)b_{\bbeta_1}(s,X(s))}{\sigma^2(s,X(s))}dX(s)
-\frac{1}{2}\int_{t-h}^t\frac{\phi_{\bxi_1}^2(s)b_{\bbeta_1}^2(s,X(s))}{\sigma^2(s,X(s))}ds\right)
\bigg |\mathcal F_{t-h}\right]
\label{eq:postpred}
\end{align}
be the posterior predictive density.
Further, for any Borel set $A$ such that $\pi(A)>0$, let
\begin{align}
&\hat f_{t-h,t,A}\left(X_{t-h,t}\right)\notag\\
&=\int_{\btheta_1\in A}
\exp\left(\int_{t-h}^t\frac{\phi_{\bxi_1}(s)b_{\bbeta_1}(s,X(s))}{\sigma^2(s,X(s))}dX(s)
-\frac{1}{2}\int_{t-h}^t\frac{\phi_{\bxi_1}^2(s)b_{\bbeta_1}^2(s,X(s))}{\sigma^2(s,X(s))}ds\right)
\pi_{t-h,A}(d\btheta_1)\notag\\
&=E\left[\exp\left(\int_{t-h}^t\frac{\phi_{\bxi_1}(s)b_{\bbeta_1}(s,X(s))}{\sigma^2(s,X(s))}dX(s)
-\frac{1}{2}\int_{t-h}^t\frac{\phi_{\bxi_1}^2(s)b_{\bbeta_1}^2(s,X(s))}{\sigma^2(s,X(s))}ds\right)
\bigg |\mathcal F_{t-h},A\right],
\end{align}
where 
\[
\pi_{t,A}(d\btheta_1)=\frac{{\bI}_A(\btheta_1)\pi_t(d\btheta_1)}{\int_A \pi_t(d\btheta_1)}
\]
is the posterior restricted to the set $A$.
We assume the following:
\begin{itemize}
\item[(H8$^{\prime\prime}$)]
\begin{equation}
\underset{t}{\lim\inf}~E_{\btheta_0}\left[\tilde{\mathcal K}'_t\left(f_{\btheta_0},\hat f_{A_t(\delta)}\right)\right]\geq\delta,
\label{eq:Q_star_t}
\end{equation}
whenever 
\begin{equation}
A_t(\delta)=\left\{\btheta_1\in\bTheta:\tilde{\mathcal K}'_t(f_{\btheta_0}, f_{\btheta_1})\geq\delta\right\}.
\label{eq:A_t}
\end{equation}
\end{itemize}
We refer to (H8$^{\prime\prime}$) as the $Q^*$ property.

\subsection{Main result on convergence of Bayes factor when two individual $SDE$'s are compared}
\label{subsec:main_result}

Let $I_0\equiv 1$ and for $t> 0$, let us define, analogous to (\ref{eq:bf_T}),
\begin{equation}
I_t=\int \frac {f_{\btheta_1,t}(X_{0,t})}{f_{\btheta_0,t}(X_{0,t})}\pi (d\btheta_1).
\label{eq:bf_t}
\end{equation}
The following lemma, proved in section S-4 of the supplement, will prove useful in proving our main theorem on convergence of Bayes factor.
\begin{lemma}
\label{lemma:towards_martingale}
\begin{align}
E_{\btheta_0}\left[\log \frac{I_t}{I_{t-h}}\bigg |\mathcal F_{t-h}\right]
&=E_{\btheta_0}\left[\log \frac{\hat f_{t-h,t}(X_{t-h,t})}{f_{\btheta_0,t-h,t}(X_{t-h,t})}\bigg |\mathcal F_{t-h}\right]
= -\tilde{\mathcal K}(f_{\btheta_0,t-h,t},\hat f_{t-h,t}).
\label{eq:bf_ratio2}
\end{align}
\end{lemma}

We make the following further assumption:
\begin{itemize}
\item[(H9$^{\prime\prime}$)]
For any $t\geq 0$, $\tilde{\mathcal K}\left(f_{\btheta_0,t,t+h_n},\hat f_{t,t+h_n}\right)$
converges in expectation 
for all sequences $\{h_n\}$ converging to zero as $n\rightarrow\infty$, with limit independent
of $\{h_n\}$. We refer to the limiting process as $\tilde{\mathcal K}'_t$. In other words, 
\begin{equation}
\underset{n\rightarrow\infty}{\lim}~
E\left(\frac{\tilde{\mathcal K}\left(f_{\btheta_0,t,t+h_n},\hat f_{t,t+h_n}\right)}{h_n}\right)
=E\left(\tilde{\mathcal K}'_t\left(f_{\btheta_0},\hat f\right)\right),
\label{eq:msd2}
\end{equation}
for any sequence $\{h_n\}$ such that $h_n\rightarrow 0$ as $n\rightarrow\infty$.
\end{itemize}
Because of Lemma \ref{lemma:towards_martingale} it follows from (H9$^{\prime\prime}$), using uniform integrability 
(which is easily seen to hold because of (H1$^{\prime\prime}$) -- (H4$^{\prime\prime}$) and (\ref{eq:moment1_T})), that
$J_{h_n}(t)=\frac{\log I_{t+h_n}-\log I_t}{h_n}$ converges in expectation  
for all sequences $\{h_n\}$ converging to zero as $n\rightarrow\infty$, with limit independent
of $\{h_n\}$. We refer to the limiting process as $J'_t$. That is, 
for any $t\geq 0$, 
\begin{equation}
\underset{n\rightarrow\infty}{\lim}~E\left(J_{h_n}(t)\right)=E\left(J'_t\right)=\frac{d}{dt}E(\log I_t).
\label{eq:msd}
\end{equation}
Now,
\begin{equation}
\hat f_{t,t+h}
=E\left[\exp\left(\int_{t}^{t+h}\frac{\phi_{\bxi_1}(s)b_{\bbeta_1}(s,X(s))}{\sigma^2(s,X(s))}dX(s)
-\frac{1}{2}\int_{t}^{t+h}\frac{\phi_{\bxi_1}^2(s)b_{\bbeta_1}^2(s,X(s))}{\sigma^2(s,X(s))}ds\right)
\bigg |\mathcal F_{t}\right],
\label{eq:postpred2}
\end{equation}
so that Lemma \ref{lemma:towards_martingale} implies 
\begin{align}
E\left[\log \frac{I_{t+h}}{I_{t}}\bigg |\mathcal F_{t}\right]
&=E\left[\log \frac{\hat f_{t,t+h}(X_{t,t+h})}{f_{\btheta_0,t,t+h}(X_{t,t+h})}\bigg |\mathcal F_{t}\right]
= -\tilde{\mathcal K}(f_{\btheta_0,t,t+h},\hat f_{t,t+h}).
\label{eq:bf_ratio3}
\end{align}
It follows, using (H9$^{\prime\prime}$), 
that
\begin{equation}
E\left(J'_t|\mathcal F_{t}\right)=-\tilde{\mathcal K}'_t\left(f_{\btheta_0},\hat f\right).
\label{eq:towards_martingale2}
\end{equation}
Note that for all sequences $\{h_n\}$ such that $h_n\rightarrow 0$ as $n\rightarrow\infty$, 
$J_{h_n}(t)=\frac{\log I_{t+h_n}-\log I_t}{h_n}$ is measurable with respect to 
$\mathcal F_{t^*}=\sigma\left(X(s):0\leq s\leq t^*\right)$, for all $t^*>t\geq 0$.
Hence, $E\left(J'_t|\mathcal F_{t}\right)\neq J'_t$.

Regarding convergence of $I_T$, we are now ready to present our main theorem whose proof is provided in section S-5 of the supplement.
\begin{theorem}
\label{theorem:bf_convergence}
Assume the $SDE$ set-up and conditions (H1$^{\prime\prime}$) -- (H9$^{\prime\prime}$).
Then 
\begin{equation}
\frac{1}{T}E_{\btheta_0}\left(\log I_T\right)\rightarrow -\delta, 
\label{eq:bf_convergence}
\end{equation}
but
\begin{equation}
\frac{1}{T^2}Var_{\btheta_0}\left(\log I_T\right)=O(1), 
\label{eq:bf_var}
\end{equation}
as $T\rightarrow\infty$.
\end{theorem}

\begin{corollary}
\label{corollary:bf_t}
For $j=1,2$, let $R_{jT}(\btheta_j)=\frac{f_{\btheta_j,T}(X_t)}{f_{\btheta_0,T}(X_T)}$, where
$\btheta_1$ and $\btheta_2$ are two different finite sets of parameters, perhaps with different dimensionalities,
associated with the two models
to be compared. For $j=1,2$, let
\[
I_{jT}=\int R_{jT}(\btheta_j)\pi_j(d\btheta_j),
\]
where $\pi_j$ is the prior on $\btheta_j$.
Let $B_T=I_{1T}/I_{2T}$ denote
the Bayes factor for comparing the two models associated with $\pi_1$ and $\pi_2$.   
Assume 
that both
the models satisfy (H1$^{\prime\prime}$) -- (H9$^{\prime\prime}$), and
have the pseudo Kullback-Leibler property 
with $\delta=\delta_1$ and $\delta=\delta_2$ respectively.
Then
\begin{equation}
E_{\btheta_0}\left(\frac{1}{T}\log B_T\right)\rightarrow \delta_2-\delta_1,
\label{eq:bf_convergence2}
\end{equation}
as $T\rightarrow\infty$.
\end{corollary}

\section{Illustration of our asymptotic result for comparing two individual $SDE$'s with a special case}
\label{sec:illustration}

Let the parameter space $\bTheta$ be compact, so that (H1$^{\prime\prime}$) holds.
Let $b_{\bbeta_j}$ and $\sigma$ satisfy
(H2$^{\prime\prime}$) such that
\begin{align}
\frac{b_{\bbeta_j}(s,x)}{\sigma(s,x)}&\equiv\eta_j(\bbeta_j);~j=0,1,
\label{eq:b_sigma_ratio}
\end{align}
so that
\begin{equation}
\frac{b_{\bbeta_1}(s,x)b_{\bbeta_0}(s,x)}{\sigma^2(s,x)}\equiv\eta_0(\bbeta_0)\eta_1(\bbeta_1).
\label{eq:b_sigma_ratio3}
\end{equation}
In the above, $\eta_1(\bbeta_1)$ 
is continuous in $\bbeta_1$.
Hence, (H3$^{\prime\prime}$) and (H4$^{\prime\prime}$) are satisfied. 
We assume that the relevant covariates and the functions $g_l$ are such that (H5$^{\prime\prime}$) holds.

Letting $\kappa_j(\bbeta_j)=\left\{\eta_j(\bbeta_j)\right\}^2$ and
$\bar{\kappa}(\bbeta_0,\bbeta_1)=\eta(\bbeta_0)\eta(\bbeta_1)$, equations 
(\ref{eq:b_sigma_ratio}) -- (\ref{eq:b_sigma_ratio3}) entail
\begin{align}
V_{\btheta_j,t}&=\kappa_j(\bbeta_j)\int_0^t\phi^2_{\bxi_j}(s)ds;\label{eq:V_example}\\
V_{\btheta_0,\btheta_1,t}&=\bar{\kappa}(\bbeta_0,\bbeta_1)\int_0^t\phi_{\bxi_0}(s)\phi_{\bxi_1}(s)ds;
\label{eq:V_example2}\\
U_{\btheta_0,t}&=\kappa_0(\bbeta_0)\int_0^t\phi^2_{\bxi_0}(s)ds
+\eta_0(\bbeta_0)\int_0^t\phi_{\bxi_0}(s)dW(s);\label{eq:U_example}\\
U_{\btheta_0,\btheta_1,t}&=\bar{\kappa}(\bbeta_0,\bbeta_1)\int_0^t\phi_{\bxi_0}(s)\phi_{\bxi_1}(s)ds
+\eta_1(\bbeta_1)\int_0^t\phi_{\bxi_1}(s)dW(s);\label{eq:U_example2}\\
V_{\btheta_j,t,t+h}&=\kappa_j(\bbeta_j)\int_t^{t+h}\phi^2_{\bxi_j}(s)ds;\label{eq:V_example_h}\\
V_{\btheta_0,\btheta_1,t,t+h}&=\bar{\kappa}(\bbeta_0,\bbeta_1)\int_t^{t+h}\phi_{\bxi_0}(s)\phi_{\bxi_1}(s)ds;
\label{eq:V_example_h2}\\
U_{\btheta_0,t,t+h}&=\kappa_0(\bbeta_0)\int_t^{t+h}\phi^2_{\bxi_0}(s)ds
+\eta_0(\bbeta_0)\int_t^{t+h}\phi_{\bxi_0}(s)dW(s);
\label{eq:U_example_h}\\
U_{\btheta_0,\btheta_1,t,t+h}&=\bar{\kappa}(\bbeta_0,\bbeta_1)\int_t^{t+h}\phi_{\bxi_0}(s)\phi_{\bxi_1}(s)ds
+\eta_1(\bbeta_1)\int_t^{t+h}\phi_{\bxi_1}(s)dW(s).
\label{eq:U_example_h2}
\end{align}
Due to (\ref{eq:V_example}) and (\ref{eq:V_example2}), we obtain
\begin{align}
\tilde{\mathcal K}'_t(f_{\btheta_0},f_{\btheta_1})&=\lim_{h\rightarrow 0} 
\frac{\tilde {\mathcal K}(f_{\btheta_0,t,t+h},f_{\btheta_1,t,t+h})}{h}\notag\\
&=\frac{1}{2}\frac{d}{dt}E_{\theta_0}(V_{\btheta_0,t})
-\frac{d}{dt}E_{\theta_0}(V_{\btheta_0,\btheta_1,t})
+\frac{1}{2}\frac{d}{dt}E_{\theta_0}(V_{\btheta_1,t})\notag\\
&=\frac{\phi_{\bxi_0}^2(t)}{2}\kappa_0(\bbeta_0)
-\phi_{\bxi_1}(t)\phi_{\bxi_0}(t)\bar{\kappa}(\bbeta_0,\bbeta_1)
+\frac{\phi_{\bxi_1}^2(t)}{2}\kappa_1(\bbeta_1)\notag\\
&=\frac{1}{2}\left(\phi_{\bxi_0}(t)\eta_0(\bbeta_0)-\phi_{\bxi_1}(t)\eta_1(\bbeta_1)\right)^2.
\label{eq:kl4}
\end{align}
Note that 
\begin{align}
\underset{t}{\inf}~\tilde{\mathcal K}'_t(f_{\btheta_0},f_{\btheta_1})
&\geq \frac{1}{2}~\underset{\bbeta_1\in\mathfrak B,\bxi_1\in\bGamma,\bz\in\mathfrak Z}{\inf}
~\left(\phi_{\bxi_0}(\bz)\eta_0(\bbeta_0)
-\phi_{\bxi_1}(\bz)\eta_1(\bbeta_1)\right)^2\notag\\
&=\delta.
\label{eq:kl_inf2}
\end{align}
It follows from (\ref{eq:kl_inf2}) that (H6$^{\prime\prime}$) holds for any prior on $\btheta_1$. 

Also, it follows directly from (\ref{eq:kl4}), that
\begin{equation}
\bar{\mathcal K}^{\infty}(f_{\btheta_0},f_{\btheta_1})
=\frac{\bar\phi^{(2)}_{\bxi_0}}{2}\kappa_0(\bbeta_0)-\bar\phi^{(2)}_{\bxi_0,\bxi_1}\bar{\kappa}(\bbeta_0,\bbeta_1)
+\frac{\bar\phi^{(2)}_{\bxi_1}}{2}\kappa_1(\bbeta_1),
\label{eq:kl_inf3}
\end{equation}
which is continuous in $(\bbeta_1,\bxi_1)$, due to the continuity assumption of $\eta_1(\bbeta_1)$ in $\bbeta_1$ and 
Lemma \ref{lemma:uniform_convergence}, which guarantees continuity of $\bar\phi^{(2)}_{\bxi_0,\bxi_1}$
and $\bar\phi^{(2)}_{\bxi_1}$ in $\bxi_1$.
%
%
Since the right-most side of (\ref{eq:kl_inf3}) is a continuous function of $\btheta_1$, 
it follows that (H7$^{\prime\prime}$) is clearly satisfied if the prior $\pi$ is dominated by the Lebesgue measure.

We now verify the $Q^*$ property (H8$^{\prime\prime}$). 
Recall that $\hat f_{A_t(\delta)}=\hat f_{t,t+h}$, since $\pi\left(A_t(\delta)\right)=1$. Since
\begin{equation}
\hat f_{t,t+h}(X_{t,t+h})\leq \underset{\btheta_1\in A_t(\delta)}{\sup}~f_{\btheta_1,t,t+h}(X_{t,t+h})
=f_{\hat\btheta_1(X_{t,t+h}),t,t+h}(X_{t,t+h}),
\label{eq:q_star_example}
\end{equation}
where $\hat\btheta_1(X_{t,t+h})\in A_t(\delta)$, is the maximizer of $f_{\btheta_1,t,t+h}$ in the compact set $A_t(\delta)$.
Hence,
\begin{align}
&\tilde{\mathcal K}\left(f_{\btheta_0,t,t+h},\hat f_{A_t(\delta)}\right)\notag\\
&=E_{\btheta_0}\left(\log f_{\btheta_0,t,t+h}(X_{t,t+h})\right)
-E_{\btheta_0}\left(\log \hat f_{A_t(\delta)}(X_{t,t+h})|\mathcal F_t\right)\notag\\
&\geq E_{\btheta_0}\left(\log f_{\btheta_0,t,t+h}(X_{t,t+h})\right)
-E_{\btheta_0}\left(\log f_{\hat\btheta_1(X_{t,t+h}),t,t+h}(X_{t,t+h})|\mathcal F_t\right)\notag\\
&=E_{\btheta_0}\left(\log\frac{f_{\btheta_0,t,t+h}(X_{t,t+h})}{f_{\hat\btheta_1(X_{t,t+h}),t,t+h}(X_{t,t+h})}
\bigg |\mathcal F_t\right)\notag\\
&=E_{\hat\btheta_1(X_{t,t+h})|\btheta_0}E_{X_{t,t+h}|\hat\btheta_1(X_{t,t+h}),\btheta_0}
\left(\log\frac{f_{\btheta_0,t,t+h}(X_{t,t+h})}
{f_{\hat\btheta_1(X_{t,t+h}),t,t+h}(X_{t,t+h})}\bigg |\mathcal F_t\right). 
\label{eq:example_K}
\end{align}
Now, let $f_{\btheta_0,tt+h}(X_{t,t+h}|Y)=\frac{f_{\btheta_0,tt+h}(X_{t,t+h})}{g_{\btheta_0,tt+h}(Y)}\bI_{\left\{\hat\btheta_1(X_{t,t+h})=Y\right\}}(Y)$ be the conditional
density of $X_{t,t+h}$ given $Y=\hat\btheta_1(X_{t,t+h})$, the latter having density $g_{\btheta_0,tt+h}(Y)$. The dominating probability measure associated with this conditional density
is $P_{0,t,t+h}$, which is the same dominating probability measure associated with $f_{\btheta_0,t,t+h}$. Then as in \ctn{Maitra15a} we have
\begin{align}
&E_{X_{t,t+h}|Y,\btheta_0}\left(\log\frac{f_{\btheta_0,t,t+h}(X_{t,t+h})}{f_{\hat\btheta_1(X_{t,t+h}),t,t+h}(X_{t,t+h})}\bigg |\mathcal F_t\right)\notag\\
&=\int\log\left(\frac{f_{\btheta_0,t,t+h}(X_{t,t+h}|Y)}{f_{Y,t,t+h}(X_{t,t+h})}\right)f_{\btheta_0,t,t+h}(X_{t,t+h}|Y)dP_{0,t,t+h}+\log g_{\btheta_0,tt+h}(Y).
\label{eq:kl_example1}
\end{align}
Since the first term of (\ref{eq:kl3}) is the Kullback-Leibler divergence between $f_{\btheta_0,t,t+h}(X_{t,t+h}|Y)$ and $f_{Y,t,t+h}(X_{t,t+h})$, it is positive for almost all $Y$.
Hence, 
\begin{equation}
E_{Y|\btheta_0}E_{X_{t,t+h}|Y,\btheta_0}\left(\log\frac{f_{\btheta_0,t,t+h}(X_{t,t+h}|Y)}{f_{Y,t,t+h}(X_{t,t+h})}\right)>0.
\label{eq:kl_positive1}
\end{equation}
Also, by Jensen's inequality, $E_{Y|\btheta_0}\left[\log g_{\btheta_0,tt+h}(Y)\right]\geq -\log E_{Y|\btheta_0}\left(\frac{1}{g_{\btheta_0,tt+h}(Y)}\right)$.
Assuming the distribution of $Y$ is dominated by the Lebesgue measure, we have $E_{Y|\btheta_0}\left(\frac{1}{g_{\btheta_0,tt+h}(Y)}\right)=|A_t(\delta)|=\int_{A_t(\delta)}dy$.
As in \ctn{Maitra15a},
once again we argue that the compact space $\bTheta$ can be rescaled appropriately with respect to suitable reparameterization such that for all $h>0$, 
$\underset{t}{\sup}~|A_t(\delta)|<\exp(-\delta h)$. 
Hence, $E_{Y|\btheta_0}\left[\log g_{\btheta_0,tt+h}(Y)\right]\geq \delta h$, which finally implies in accordance with (\ref{eq:kl_positive1}), that
$\tilde{\mathcal K}\left(f_{\btheta_0,t,t+h},\hat f_{A_t(\delta)}\right)\geq\delta h$. Hence, $\tilde{\mathcal K}'\left(f_{\btheta_0},\hat f_{A_t(\delta)}\right)\geq\delta$,
showing that the $Q^*$ property is satisfied.


To see that (H9$^{\prime\prime}$) holds, first observe that it follows from the proof of Lemma \ref{lemma:towards_martingale}
that 
$\frac{I_{t+h}}{I_{t}}=\frac{\hat f_{t,t+h}}{f_{\btheta_0,t,t+h}}$, which implies
\begin{align}
\frac{\tilde{\mathcal K}\left(f_{\btheta_0,t,t+h},\hat f_{t,t+h}\right)}{h}
&=\frac{E_{\btheta_0}\left(\log f_{\btheta_0,t,t+h}(X_{t,t+h})-\log \hat f_{t,t+h}(X_{t,t+h})|\mathcal F_t\right)}{h}
\label{eq:verify_H11}
\end{align}
Now,
\begin{equation}
E_{\btheta_0}\left(\log f_{\btheta_0,t,t+h}(X_{t,t+h})|\mathcal F_t\right)
=\frac{\kappa_0(\bbeta_0)}{2}\int_t^{t+h}\phi^2_{\bxi_0}(s)ds=\frac{\kappa_0(\bbeta_0)}{2}h\phi^2_{\bxi_0}(s^*(h)),
\label{eq:inner0}
\end{equation}
by the mean value theorem for integrals, where $s^*(h)\rightarrow t$, as $h\rightarrow 0$.
Hence, using continuity of $\phi_{\bxi_0}(t)$ in $t$, we obtain 
\begin{equation}
\underset{h\rightarrow 0}{\lim}~\frac{1}{h}~E_{\btheta_0}\left(\log f_{\btheta_0,t,t+h}(X_{t,t+h})|\mathcal F_t\right)
=\frac{\kappa_0(\bbeta_0)}{2}~ \underset{h\rightarrow 0}{\lim}~\phi^2_{\bxi_0}(s^*(h))
=\frac{\kappa_0(\bbeta_0)}{2}\phi^2_{\bxi_0}(t).
\label{eq:verify_H11_1}
\end{equation}
To deal with $E_{\btheta_0}\left(\log \hat f_{t,t+h}(X_{t,t+h})|\mathcal F_t\right)$, 
note that
for any $X_{t,t+h}$, by the mean value theorem for integrals, 
$$\hat f_{t,t+h}(X_{t,t+h})=f_{\breve\btheta(X_{0,t},X_{t,t+h}),t,t+h}(X_{t,t+h}),$$ where 
$\breve\btheta(X_{0,t},X_{t,t+h})\in\bTheta$. 
It is clear that
$\breve\btheta(X_{0,t},X_{t,t+h})\rightarrow \breve\btheta(X_{0,t},X_{t})=\breve\btheta(X_{0,t})$ almost surely, 
as $h\rightarrow 0$. 
Hence,
\begin{align}
&E_{\btheta_0}\left(\log \hat f_{t,t+h}(X_{t,t+h})|\mathcal F_t\right)\notag\\
&=E_{\breve\btheta(X_{0,t},X_{t,t+h})|\btheta_0}E_{X_{t,t+h}|\breve\btheta(X_{0,t},X_{t,t+h})=\alpha,\btheta_0}
\left(\log f_{\{\breve\btheta(X_{0,t},X_{t,t+h})=\alpha\},t,t+h}(X_{t,t+h})|\mathcal F_t\right),\notag
\end{align}
where
\begin{align}
&E_{\btheta_0}\left(\log \hat f_{t,t+h}(X_{t,t+h})|\mathcal F_t\right)\notag\\
&=E_{\btheta_0}\left(\log f_{\{\breve\btheta(X_{0,t},X_{t,t+h})\},t,t+h}(X_{t,t+h})|\mathcal F_t\right)\notag\\
&=E_{\btheta_0}\left(\bar\kappa(\bbeta_0,\breve\bbeta_1(X_{0,t},X_{t,t+h}))
\int_t^{t+h}\phi_{\breve\bxi_1(X_{0,t},X_{t,t+h})}(s)\phi_{\bxi_0}(s)ds\right.\notag\\
&\quad\quad\quad\left.-\frac{\kappa_1(\breve\bbeta_1(X_{0,t},X_{t,t+h}))}{2}\int_t^{t+h}
\phi^2_{\breve\bxi_1(X_{0,t},X_{t,t+h})}(s)ds\right)
\notag\\
&=E_{\btheta_0}\left(\bar\kappa(\bbeta_0,\breve\bbeta_1(X_{0,t},X_{t,t+h}))
h\phi_{\breve\bxi_1(X_{0,t},X_{t,t+h})}(s_1(h))\phi_{\bxi_0}(s_1(h))\right.\notag\\
&\quad\quad\quad\left.-\frac{\kappa_1(\breve\bbeta_1(X_{0,t},X_{t,t+h}))}{2}h\phi^2_{\breve\bxi_1(X_{0,t},X_{t,t+h})}(s_2(h))\right),
\label{eq:inner1}
\end{align}
where $s_1(h),s_2(h)\in[t,t+h]$, associated with the mean value theorem for integrals. Hence,
$s_1(h)\rightarrow t$ and $s_2(h)\rightarrow t$, almost surely as $h\rightarrow 0$.

Continuity of $\bar\kappa(\cdot),\kappa_1(\cdot)$ and the results
$\breve\btheta(X_{0,t},X_{t,t+h})\rightarrow \breve\btheta(X_{0,t})$, $s_1(h)\rightarrow t$, $s_2(h)\rightarrow t$, almost surely, 
as $h\rightarrow 0$, in conjunction with the dominated convergence theorem exploiting boundedness of the functions 
$\phi_{\bxi_0}$, $\phi_{\breve\bxi_1}$, $\bar\kappa$ and $\kappa_j$, imply, using continuity of 
$\phi_{\breve\bxi_1}(t)$ in $t$, that
\begin{align}
&\underset{h\rightarrow 0}{\lim}~\frac{1}{h}~E_{\btheta_0}\left(\log \hat f_{t,t+h}(X_{t,t+h})|\mathcal F_t\right)\notag\\ 
&\quad\quad=\phi_{\bxi_0}(t)\phi_{\breve\bxi_1(X_{0,t})}(t)\bar\kappa(\bbeta_0,\breve\bbeta(X_{0,t}))
-\frac{\phi_{\breve\bxi_1(X_{0,t})}^2(t)}{2}\kappa_1(\breve\bbeta(X_{0,t})).\notag
\end{align}
In other words, the limit of (\ref{eq:verify_H11}) exists and is unique as $h\rightarrow 0$.
Now, equations (\ref{eq:inner0}) and (\ref{eq:inner1}) along with dominated convergence theorem imply that  (H9$^{\prime\prime}$) holds. 

Thus, all the assumptions required for Theorem \ref{theorem:bf_convergence} and Corollary \ref{corollary:bf_t}
are satisfied.
Hence, both (\ref{eq:bf_convergence}) and (\ref{eq:bf_convergence2}) hold.

\section{Asymptotic convergence of Bayes factor in the $SDE$ set-up with respect to  
number of individuals and time}
\label{sec:asymp_BF_nT}


\subsection{Convergence of Bayes factor in the $iid$ set-up}
\label{subsec:iid}


Although Theorem \ref{theorem:bf_convergence} fails to ensure consistency of the Bayes factor as $T\rightarrow\infty$
in the sense that the relevant variance is asymptotically positive,
the theorem is useful to prove almost sure consistency when $T\rightarrow\infty$ as well as $n\rightarrow\infty$, for
both $iid$ and non-$iid$ situations. Theorem \ref{theorem:bf_convergence_iid} formalizes this for the $iid$ set-up,
while Theorem \ref{theorem:bf_convergence_non_iid} establishes almost sure consistency of the Bayes factor
in the non-$iid$ situation. Proofs of these theorems are contained in section S-6 and S-9 respectively in the supplement.
\begin{theorem}
\label{theorem:bf_convergence_iid}
Assume the $iid$ set-up;
also assume that conditions (H1$^{\prime\prime}$) -- (H9$^{\prime\prime}$) hold for each $SDE$ 
in the systems (\ref{eq:sde3}) and (\ref{eq:sde4}). Then
\begin{equation}
\frac{1}{nT}\log\tilde I_{n,T}
\rightarrow -\delta,
\label{eq:bf_convergence_iid2}
\end{equation}
almost surely, as $n\rightarrow\infty$ and $T\rightarrow\infty$.
\end{theorem}

The following corollary is obvious.
\begin{corollary}
\label{corollary:bf_t_iid}
For $j=1,2$, and $i=1,\ldots,n$, 
let $R_{j,i,T}(\btheta^{(i)}_{j})=\frac{f_{\btheta^{(i)}_j,i,T}(X_{i,0,T})}{f_{\btheta^{(i)}_0,i,T}(X_{i,0,T})}$, where, for each $i$,
$\btheta^{(i)}_1$ and $\btheta^{(i)}_2$ are two different finite sets of parameters, perhaps with different dimensionalities,
associated with the two systems (\ref{eq:sde3}) and (\ref{eq:sde4})
to be compared. For $j=1,2$, let
\[
\tilde I_{j,n,T}=\prod_{i=1}^n\int R_{j,i,T}(\btheta^{(i)}_j)\pi_j(d\btheta^{(i)}_j),
\]
where $\pi_j$ is the prior on $\btheta^{(i)}_j$, for $i=1,2,\ldots$.
Let $B_{n,T}=\tilde I_{1,n,T}/\tilde I_{2,n,T}$ denote
the Bayes factor for comparing the two models associated with $\pi_1$ and $\pi_2$.   
Assume the $iid$ case and suppose 
that both
the systems satisfy (H1$^{\prime\prime}$) -- (H9$^{\prime\prime}$), and
have the pseudo Kullback-Leibler property 
with $\delta=\delta_1$ and $\delta=\delta_2$ respectively.
Then
\begin{equation*}
\frac{1}{nT}\log B_{n,T}\rightarrow \delta_2-\delta_1,
\label{eq:bf_convergence2_iid}
\end{equation*}
almost surely, as $n\rightarrow\infty$ and $T\rightarrow\infty$.
\end{corollary}

\subsection{Convergence of Bayes factor in the non-$iid$ set-up}
\label{subsec:bf_convergence_non_iid}

We now relax the assumptions $x^i=x$ and $\xi^{(i)}_{1j}=\xi^{(i)}_{2j}=\xi^{(i)}_{3j}=\cdots=\xi^{(i)}_{pj}=0$ for $j=0,1$. 
Thus, we are now in a non-$iid$ situation
where the processes $X_i(\cdot);~i=1,\ldots,n$, are independently,
but not identically distributed.
As mentioned in Section \ref{subsec:covariates} we assume that $\btheta^{(i)}_1\stackrel{iid}{\sim}\pi$.
In this set-up, for each $\bz\in\mathfrak Z=\left\{\bz(t)\in\boldcal Z:t\in[0,\infty)\right\}$,
it holds, due to Theorem \ref{theorem:bf_convergence}, that
\begin{equation}
\frac{1}{T}E_{\btheta_0}\left(\log I_{x,T,\bz}\right)\rightarrow -\delta(x,\bz), 
\label{eq:bf_convergence_non_iid}
\end{equation}
as $T\rightarrow\infty$, where $\delta(x,\bz)$ depends upon the initial value $x\in\mathfrak X$ 
and the set of time-dependent covariates 
$\bz\in\mathfrak Z$. 
The following lemma shows that $\delta(x,\bz)$ is continuous in $(x,\bz)\in\mathfrak X\times\mathfrak Z$.
\begin{lemma}
\label{lemma:delta_continuity}
Assume the conditions of Theorem \ref{theorem:bf_convergence}. Then,
$\delta(x,\bz)$ is continuous in $(x,\bz)\in\mathfrak X\times\mathfrak Z$.
\end{lemma}

Now consider the following limit:
\begin{equation}
\delta^{\infty}=\underset{n\rightarrow\infty}{\lim}~\frac{1}{n}\sum_{i=1}^n\delta(x^i,\bz_i).
\label{eq:delta_convergence}
\end{equation}
The following lemma shows that the above limit exists for all sequences 
$\left\{(x^i,\bz_i)\right\}_{i=1}^{\infty}\in\mathfrak X\times\mathfrak Z$.
\begin{lemma}
\label{lemma:delta_limit_existence}
The limit (\ref{eq:delta_convergence}) exists for all sequences $\left\{(x^i,\bz_i)\right\}_{i=1}^{\infty}
\in\mathfrak X\times\mathfrak Z$.
\end{lemma}

Proof of these two lemmas are provided in section S-7 and S-8 respectively in the supplement.
Now, we  have the following theorem.
\begin{theorem}
\label{theorem:bf_convergence_non_iid}
Assume the non-$iid$ set-up, and conditions (H1$^{\prime\prime}$) -- (H9$^{\prime\prime}$), for each $SDE$ 
in the systems (\ref{eq:sde3}) and (\ref{eq:sde4}). Then
\begin{equation}
\frac{1}{nT}\log\tilde I_{n,T}
\rightarrow -\delta^{\infty},
\label{eq:bf_convergence_non_iid3}
\end{equation}
almost surely, as $T\rightarrow\infty$ and $n\rightarrow\infty$.
\end{theorem}

We then have the following corollary for the non-$iid$ case.
\begin{corollary}
\label{corollary:bf_t_non_iid}
For $j=1,2$, and $i=1,\ldots,n$, 
let $R_{j,i,T}(\btheta^{(i)}_{j})=\frac{f_{\btheta^{(i)}_j,i,T}(X_{i,0,T})}{f_{\btheta^{(i)}_0,i,T}(X_{i,0,T})}$, 
where, for each $i$,
$\btheta^{(i)}_1$ and $\btheta^{(i)}_2$ are two different finite sets of parameters, perhaps with different dimensionalities,
associated with the two systems (\ref{eq:sde3}) and (\ref{eq:sde4})
to be compared. For $j=1,2$, let
\[
\tilde I_{j,n,T}=\prod_{i=1}^n\int R_{j,i,T}(\btheta^{(i)}_j)\pi_j(d\btheta^{(i)}_j),
\]
where $\pi_j$ is the prior on $\btheta^{(i)}_j;~i=1,2,\ldots$.
Let $B_{n,T}=\tilde I_{1,n,T}/\tilde I_{2,n,T}$ denote
the Bayes factor for comparing the two models associated with $\pi_1$ and $\pi_2$.   
Assume the non-$iid$ case and suppose 
that both
the systems satisfy (H1$^{\prime\prime}$) -- (H9$^{\prime\prime}$), and
have the pseudo Kullback-Leibler property 
with $\delta_i=\delta_{1i}$ and $\delta_i=\delta_{2i}$ respectively.
Let, for $j=1,2$,
\begin{equation*}
\delta^{\infty}_j=\underset{n\rightarrow\infty}{\lim}~\frac{1}{n}\sum_{i=1}^n\delta_{ji}.
\label{eq:delta_convergence2}
\end{equation*}
Then
\begin{equation*}
\frac{1}{nT}\log B_{n,T}\rightarrow \delta^{\infty}_2-\delta^{\infty}_1,
\label{eq:bf_convergence2_non_iid}
\end{equation*}
almost surely, as $n\rightarrow\infty$ and $T\rightarrow\infty$. 
\end{corollary}

\section{Simulation studies}
\label{sec:simulated_data}

\subsection{Covariate selection when $n=1$, $T=5$}
\label{subsec:N1T5}

We first demonstrate with simulation study the finite sample analogue of Bayes factor analysis
associated with a single individual, when $T\rightarrow\infty$ . In this regard, 
we consider modeling a single individual by
\begin{equation}
dX(t)=(\xi_1+\xi_2z_1(t)+\xi_3z_2(t)+\xi_4z_3(t))(\xi_5+\xi_6X(t))dt+\sigma dW(t),
\label{eq:sde2_appl}
\end{equation}
where we fix our diffusion coefficient as $\sigma=20$. 
We consider the initial value $X(0)=0$ and the time interval $[0, T]$ with $T=5$. 

To achieve numerical stability of the marginal likelihood corresponding to data we choose 
the true values of $\xi_i$; $i=1,\ldots,6$ as follows:
$\xi_{i}\stackrel{iid}{\sim} N(\mu_i,0.001^2)$, where 
$\mu_i\stackrel{iid}{\sim} N(0,1)$. This is not to be interpreted as the prior; this
is just a means to set the true values of the parameters of the data-generating model.
\\[2mm]
We assume that the time dependent covariates $z_i(t)$ satisfy the following $SDE$s 
\begin{align}
dz_1(t)=&(\tilde\theta_{1}+\tilde\theta_{2}z_1(t))dt+ dW_1(t)\notag\\
dz_2(t)=&\tilde\theta_{3}dt+ dW_2(t)\notag\\
dz_3(t)=&\tilde\theta_{4}z_3(t))dt+ dW_3(t),
\label{eq:covariate_appl}
\end{align}
where $W_i(\cdot)$; $i=1,2,3$, are independent Wiener processes, and 
$\tilde\theta_{i}\stackrel{iid}{\sim} N(0,0.01^2)$ for $i=1,\cdots,4$. 

We obtain the covariates by first simulating $\tilde\theta_{i}\stackrel{iid}{\sim} N(0,0.01^2)$ for $i=1,\cdots,4$,
fixing the values, and then by simulating the covariates using the $SDE$s (\ref{eq:covariate_appl}) 
by discretizing the time interval $[0,5]$ into $500$ equispaced time points. 
In all our applications we have standardized the covariates over time so that they have zero means
and unit variances.

Once the covariates are thus obtained, we assume that the data are generated from the (true) model where all the covariates
are present. For the true values of the parameters, we simulated $(\xi_1,\ldots,\xi_6)$ from
the prior and treated the obtained values as the true set of parameters $\btheta_0$.
We then generated the data using (\ref{eq:sde2_appl}) by discretizing the time interval $[0,5]$ 
into $500$ equispaced time points.
%

As we have three covariates so we will have $2^3=8$ different models. Denoting
a model by the presence and absence of the respective covariates, it then is the case that 
$(1,1,1)$ is the true, data-generating model, while $(0,0,0)$, $(0,0,1)$, $(0,1,0)$, $(0,1,1)$, $(1,0,0)$,
$(1,0,1)$, and $(1,1,0)$ are the other $7$ possible models. 

As per our theory, for a single individual, the Bayes factor is not consistent for increasing time domain.
However, we have shown that 
\begin{equation*}
\frac{1}{T}E_{\btheta_0}(\log I_T)\rightarrow -\delta
\label{eq:bf_appl}
\end{equation*}
as $T\rightarrow \infty$. Thus, the Bayes factor is consistent with respect to the expectation. 
Our simulation results show that this holds even for the time domain $[0,5]$,
where we approximate the expectation with the average of 
$1000$ realizations of $I_T$
associated with as many simulated data sets. 

\subsubsection{Case 1: the true parameter set $\btheta_0$ is fixed}
\label{subsubsec:case1_appl2}
{\bf Prior on $\btheta$}
\\[2mm]
We first obtain the maximum likelihood estimator ($MLE$) of $\btheta$ using simulated annealing
and then consider a normal prior with the $MLE$ as the mean and variance $0.8^2\mathbb I_6$, where
$\mathbb I_6$ is the identity matrix of order $6$.
\\[2mm]
{\bf Form of the Bayes factor}
\\[2mm]
In this case the related Bayes factor has the form
\begin{equation}
I_T=\int \frac {f_{\btheta_1,T}(X_{0,T})}{f_{\btheta_0,T}(X_{0,T})}\pi (d\btheta_1),
\label{eq:bf_T_appl2}
\end{equation}
where $\btheta_0=(\xi_{0,1},\xi_{0,2},\xi_{0,3},\xi_{0,4},\xi_{0,5},\xi_{0,6})$ is the true parameter set 
and $\btheta_1=(\xi_1,\xi_2,\xi_3,\xi_4,\xi_5,\xi_6)$ is the unknown set of parameters corresponding to 
any other model. 
Table \ref{table:values2} describes the results of our Bayes factor analyses. 
\begin{table}[h]
\centering
\caption{Bayes factor results}
\label{table:values2}
\begin{tabular}{|c||c|}
\hline
Model & Averaged $\frac{1}{5}\log I_5$\\
\hline
$(0,0,0)$ & -2.5756029\\ 
$(0,0,1)$ & -0.913546\\
$(0,1,1)$ & -0.5454860\\
$(0,1,0)$ & -0.763952\\
$(1,0,0)$ & -2.5774163\\
$(1,0,1)$ & -0.9312218\\
$(1,1,0)$ & -0.7628154\\
\hline
\end{tabular}
\end{table}
It is clear from the 7 values of the table that the correct model $(1,1,1)$ is always preferred.

\subsubsection{Case 2: the parameter set $\btheta_0$ is random and has the prior distribution $\pi$}
\label{subsubsec:case2_appl2}
As before, we consider the same form of the prior as in Section \ref{subsubsec:case1_appl2},
but with variance $0.1^2\mathbb I_6$. 
In this case we calculate marginal likelihood of the 8 possible models, and approximate
$$\frac{1}{5}E_{\btheta_0}\left(\log\int f_{i,\btheta_1,5}(X_{0,5})\pi (d\btheta_1)\right)$$ for $i=1,\ldots,8$ by 
averaging over $1000$ replications of the data obtained from the true model. 
Denoting its values by $\ell_i$,  Table \ref{table:deltavalues2} shows that $\ell_8$ 
is the highest, implying consistency of the averaged Bayes factor.
\begin{table}[h]
\centering
\caption{Averages of $\frac{1}{5}\times$ marginal log-likelihood}
\label{table:deltavalues2}
\begin{tabular}{|c||c|}
\hline
Model & $\ell_i$\\
\hline
$(0,0,0)$ & -1.21923\\ 
$(0,0,1)$ & -0.21428\\
$(0,1,0)$ & 1.47992\\
$(0,1,1)$ & 2.102966\\
$(1,0,0)$ & -1.222362\\
$(1,0,1)$ & -0.21898\\
$(1,1,0)$ & 1.459921\\
$(1,1,1)$ & 2.121237 ($\mbox{true model}$)\\
\hline
\end{tabular}
\end{table}

\subsection{Bayes factor analysis for $n=15$ and $T=5$}
\label{subsec:N15T5}

In this case we allow our parameter and the covariate sets to vary from individual to individual. 
We consider $15$ individuals modeled by
\begin{equation}
dX_i(t)=(\xi_1^i+\xi_2^iz_1(t)+\xi_3^iz_2(t)+\xi_4^iz_3(t))(\xi_5^i+\xi_6^iX_i(t))dt+\sigma_i dW_i(t)
\label{eq:sde3_appl3}
\end{equation}
for $i=1,\cdots,15$. We fix our diffusion coefficients as $\sigma_{i+1}=\sigma_i +5$ for $i=1\cdots,14$ 
where $\sigma_1=10$. We consider the initial value $X(0)=0$ and the interval $[0, T]$, with $T=5$. 
As before, we generated the observed data after discretizing the time interval into $500$ equispaced time points. 
Here our covariates and the parameter set 
$\btheta_0^i=(\xi_{0,1}^i,\xi_{0,2}^i,\xi_{0,3}^i,\xi_{0,4}^i,\xi_{0,5}^i,\xi_{0,6}^i)$; $i=1,\ldots,15$, 
are simulated in a similar way as mentioned in Section \ref{subsec:N1T5}.

For each of the $15$ individuals, the true set of covariate combination is randomly selected. 
Thus, for a given model, there are $15$ sets of covariate combinations to be compared with other models
consisting of $15$ different sets of covariate combinations.
To decrease computational burden we compare the true model with $100$ other models consisting of 
different sets of covariate combinations.  
%

The Bayes factor corresponding to the $j$-th covariate combination is given by
\begin{equation}
I_{nT}^j=\prod_{i=1}^n\int \frac{f^j_{i,\btheta^{(i)}_1}(X_{i,0,T})}
{f_{i,\btheta^{(i)}_0}(X_{i,0,T})}\pi\left(\btheta^{(i)}_{1}\right)d\btheta^{(i)}_{1}
\label{eq:I_iT_appl}
\end{equation}
for $j=1,\cdots,100$, where $n=15$, $T=5$ and $\btheta_0^{(i)}$ is the true parameter set corresponding to 
the $i$-th individual.

We obtain the $MLE$ of the $15$ parameter sets by simulated annealing. Then we calculate the Bayes factor 
with the prior such that the parameter components are independent normal with means as the respective $MLE$s and variances $1$. 
In all the cases corresponding to $100$ covariate combinations we obtain 
$\frac{1}{nT}\log I^j_{nT}<0$ for $j=1,\cdots,100$. 
Thus, Bayes factor indicated the correct covariate combination in all the cases considered. 
We also considered the case when a normal prior is considered for the parameters of the true model. 
In this case with respect to the component-wise independent normal prior with individual mean as obtained 
from simulated annealing and component-wise variance $0.1^2$, we obtain 
\begin{equation}
\frac{1}{15\times 5}\left[\log\left(\prod_{i=1}^{15}\int f^j_{i,\btheta_1^{(i)}}(X_{i,0,T})
\pi(\btheta^{(i)}_1)d\btheta_1^{(i)}\right)-\log\left(\prod_{i=1}^{15}\int f_{i,\btheta_0^{(i)}}(X_{i,0,T})
\pi(\btheta^{(i)}_0)d\btheta_0^{(i)}\right)\right]<0, 
\label{eq:INT_appl}
\end{equation}
for $j=1,\cdots,100$. Indeed, it turned out that 
$\frac{1}{15\times 5}\log\left(\prod_{i=1}^{15}\int f_{i,\btheta_0^{(i)}}(X_{i,0,T})
\pi(\btheta^{(i)}_0)d\btheta_0^{(i)}\right)=0.4865$ 
and the maximum of 
$\frac{1}{15\times 5}\log\left(\prod_{i=1}^{15}\int f ^j_{i,\btheta_1^{(i)}}(X_{i,0,T})
\pi(\btheta^{(i)})d\btheta^{(i)}\right)$ is $0.4127$.
In other words, the Bayes factor consistently selects the correct model even in this situation.

\section{Company-wise data from national stock exchange}
\label{sec:truedata} 

\begin{figure}
\begin{subfigure}
\centering
\includegraphics[height=6cm,width=5cm]{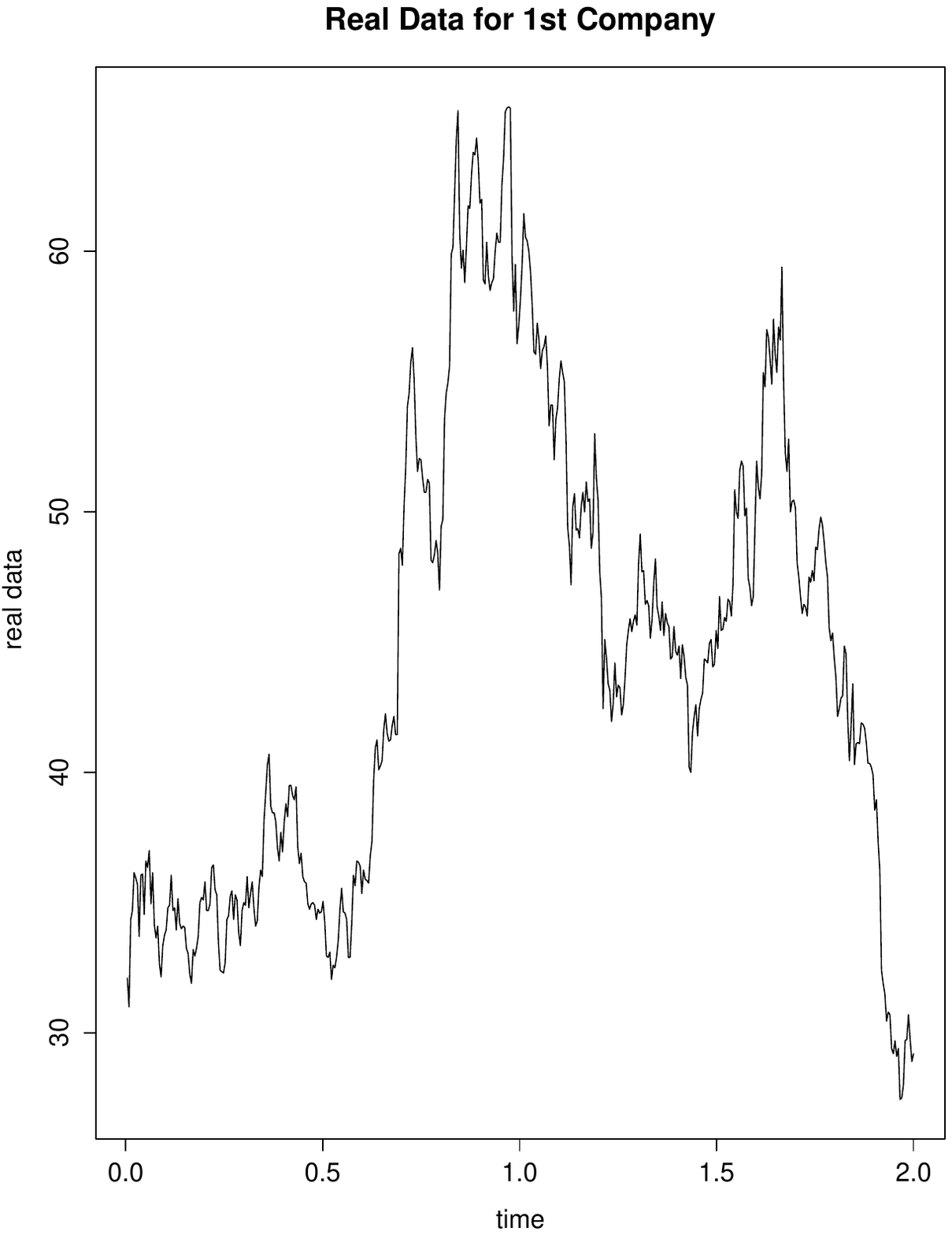}
\label{fig:sfig2}
\end{subfigure}
\begin{subfigure}
\centering
\includegraphics[height=6cm,width=5cm]{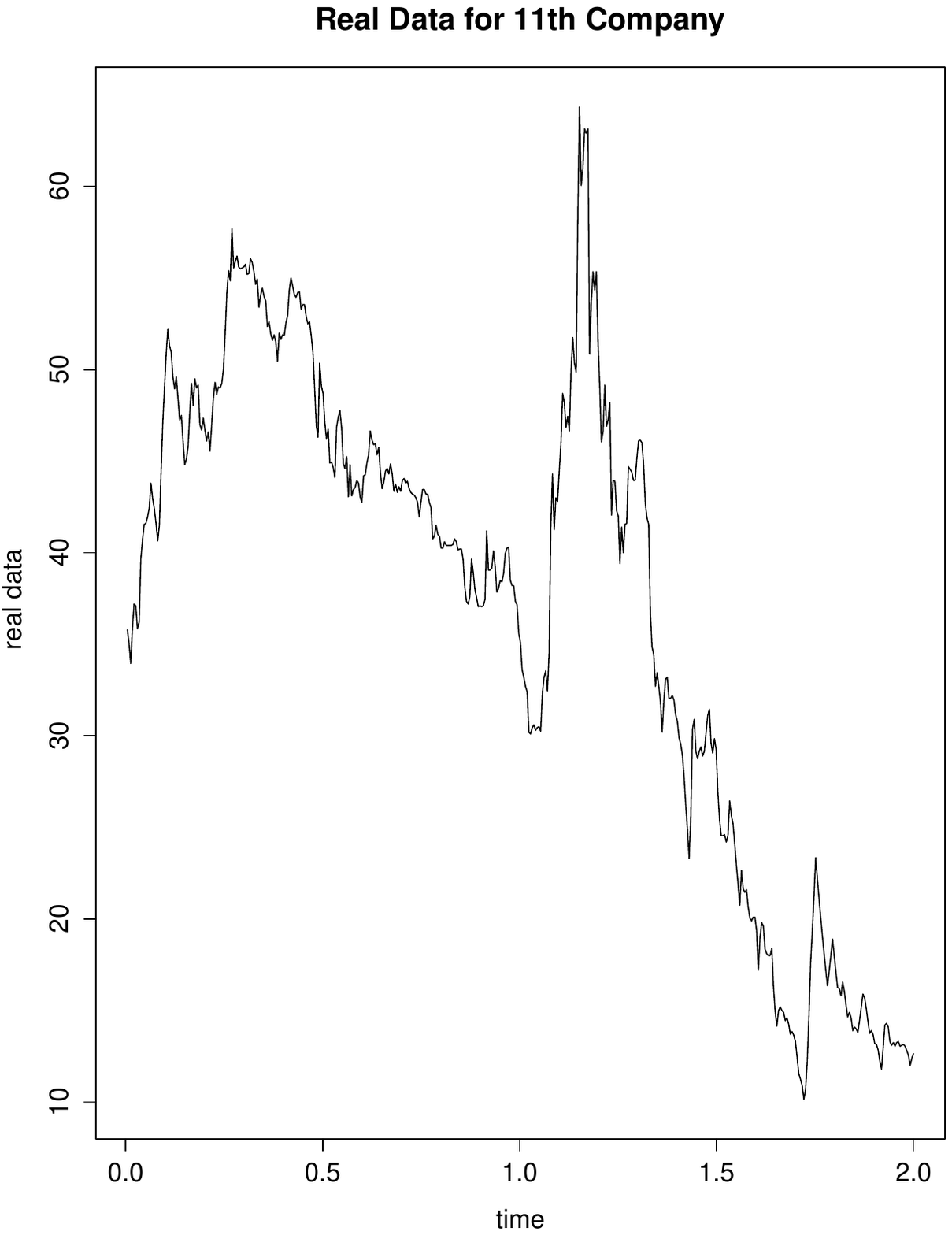}
\label{fig:sfig3}
\end{subfigure}
\begin{subfigure}
\centering
\includegraphics[height=6cm,width=5cm]{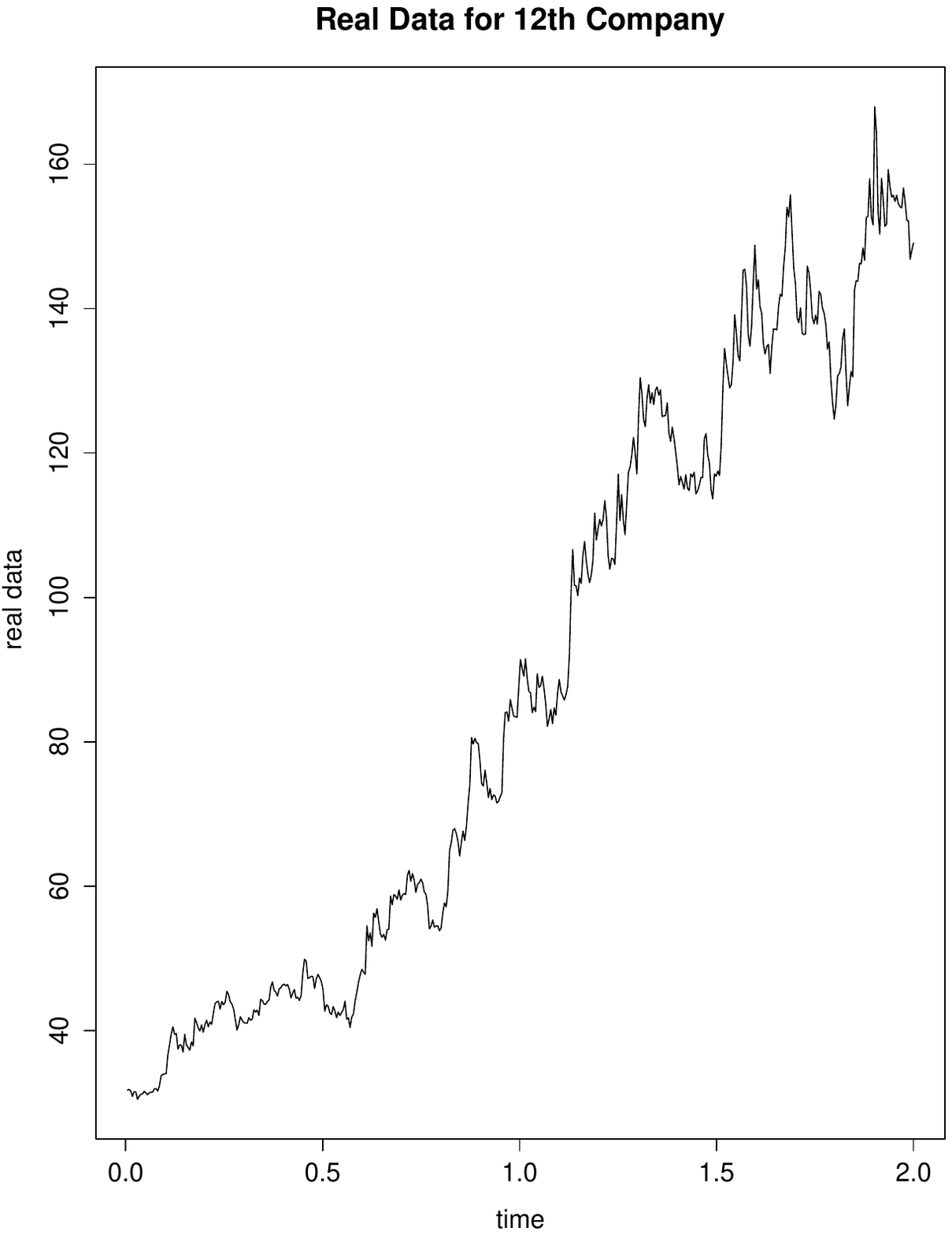}
\label{fig:sfig4}
\end{subfigure}
\caption{Some company-wise time-series data.}
\label{fig:true15}
\end{figure}

To deal with real data we collect the stock market data 
($467$ observations during the time range August $5$, 2013, to June $30$, 2015) 
for $15$ companies which is available on 
{\it www.nseindia.com}. The nature of some company-wise data are shown in Figure \ref{fig:true15}.

Each company-wise data is modeled by various availabe standard financial $SDE$ 
models with the available ``fitsde" package in $R$. After obtaining the BIC (Bayesian Information Criterion) 
for each company corresponding to each available financial model, we find that the minimum value of BIC corresponds to 
the $CKLS$ model, given, for process $X(t)$, by
$$dX(t)=(\theta_1+\theta_2X(t))dt+\theta_3X(t)^{\theta_4}dW(t).$$
As per our theory we treat the diffusion coefficient as a fixed quantity. 
So, after obtaining the estimated value of the coefficients by the ``fitsde" function, 
we fix the values of $\theta_3$ and $\theta_4$, so that the diffusion coefficient becomes fixed.
We let $\theta_3=A$, $\theta_4=B$.

In this $CKLS$ model, we now wish to include time varying covariates. 
In our work we consider the ``close price" of each company. The stock market data is assumed to 
be dependent on IIP general index, bank interest rate, US dollar exchange rate and on various other quantities. 
But we assume only these three quantities as possibly the most important time dependent covariates. 

Briefly, IIP, that is, index of industrial production, is a measurement which represents 
the status of production in the industrial sector for a given period of time compared to a reference period of time. 
It is one of the best statistical data, which helps us measure the level of industrial activity in Indian economy. 
Its importance lies in the fact that low industrial production  will result in lower corporate sales and profits, 
which will directly affect stock prices. So a direct impact of weak IIP data is a sudden fall in stock prices. 

As the IIP data is purely industrial data, banking sector is not included in it. So, 
we also consider the bank interest rate as another covariate. Note that, higher the bank interest rate, 
fixed deposits become more attractive and one will preferably deposit money in bank rather than invest in stock market. 
Besides, companies with a high amount of loans in their balance sheets would be affected very seriously. 
Interest cost on existing debt would go up affecting their EPS (Earning per Share) and ultimately the stock prices. 
But during low interest rate these companies would stand to gain. Banking sector is likely to benefit most due 
to high interest rates. The Net Interest Margins (it is the difference between the interest they 
earn on the money they lend and the interest they pay to the depositors) for banks is likely to 
increase leading to growth in profits and the stock prices. Hence, it is clear that, the interest rates 
and stock markets are inversely related. As the interest rates go up, stock market 
activities tend to come down. 

Finally, exchange rates directly affect the realized return on 
an investment portfolio with overseas holdings. If one own stock in a foreign company and the 
local currency goes up, the value of the investment also goes up. 
Foreign investment is also related very much to US dollar exchange rate.

\begin{figure}
\begin{subfigure}
\centering
\includegraphics[height=6cm,width=5cm]{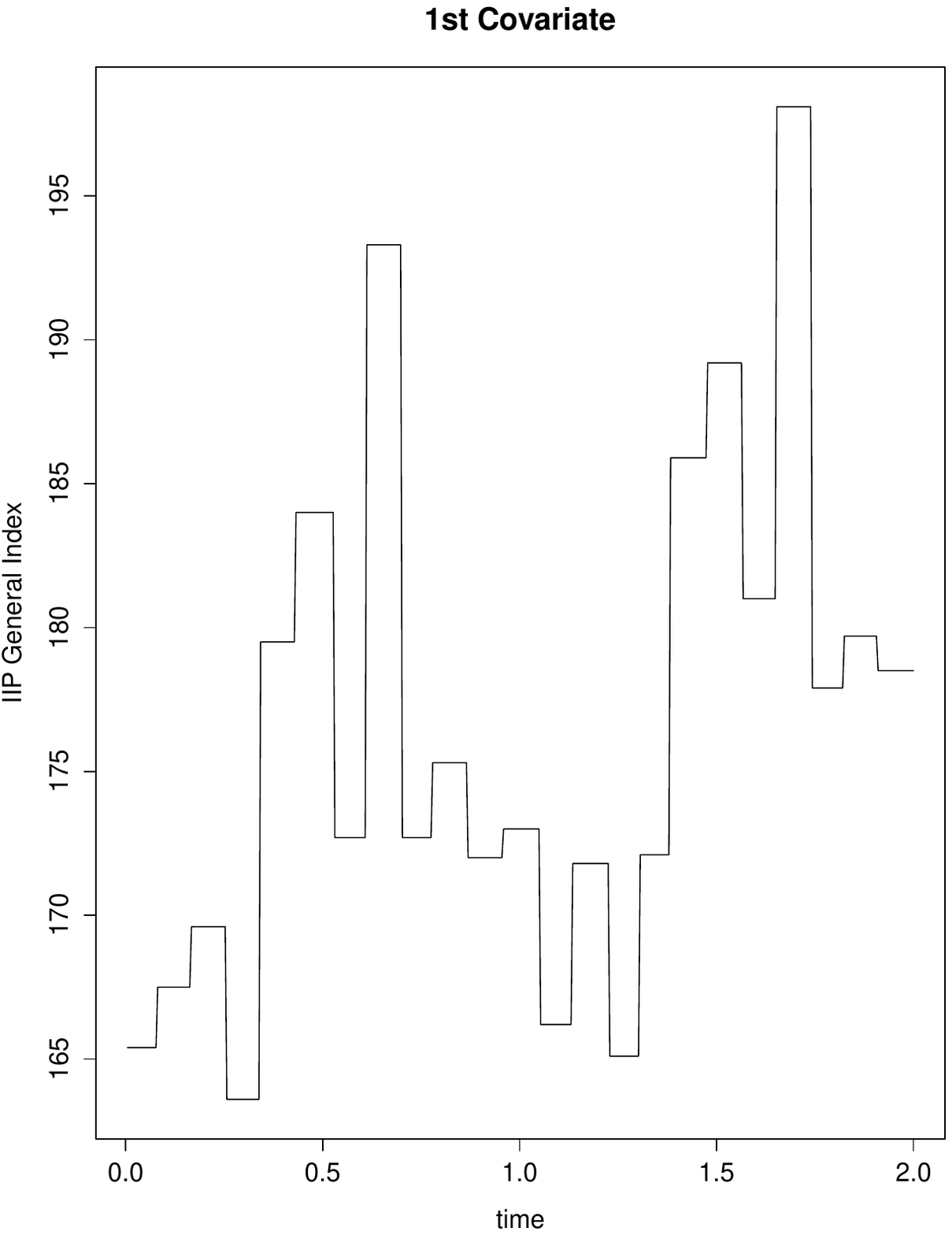}
\label{fig:realdata_sfig2}
\end{subfigure}
\begin{subfigure}
\centering
\includegraphics[height=6cm,width=5cm]{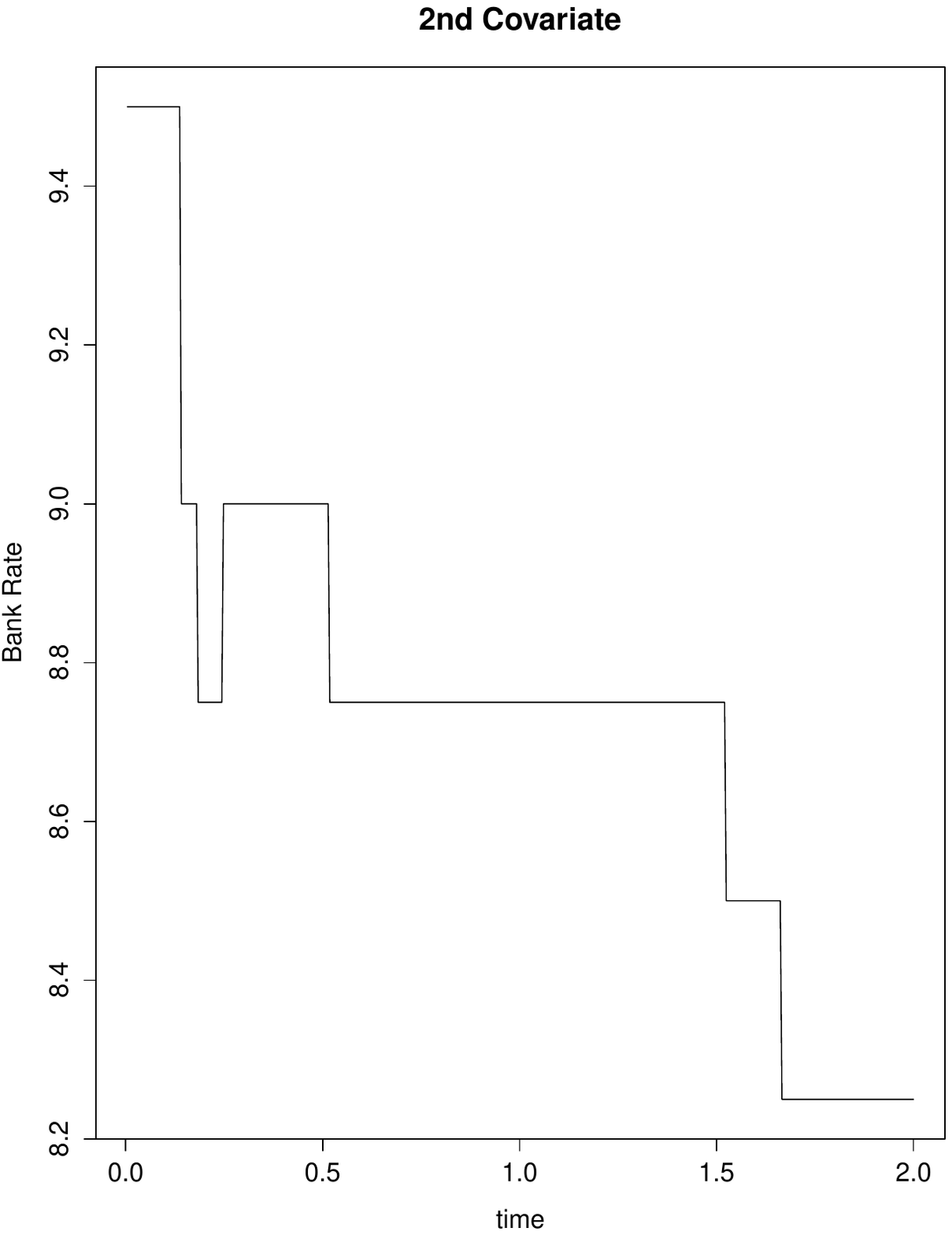}
\label{fig:realdata_sfig3}
\end{subfigure}
\begin{subfigure}
\centering
\includegraphics[height=6cm,width=5cm]{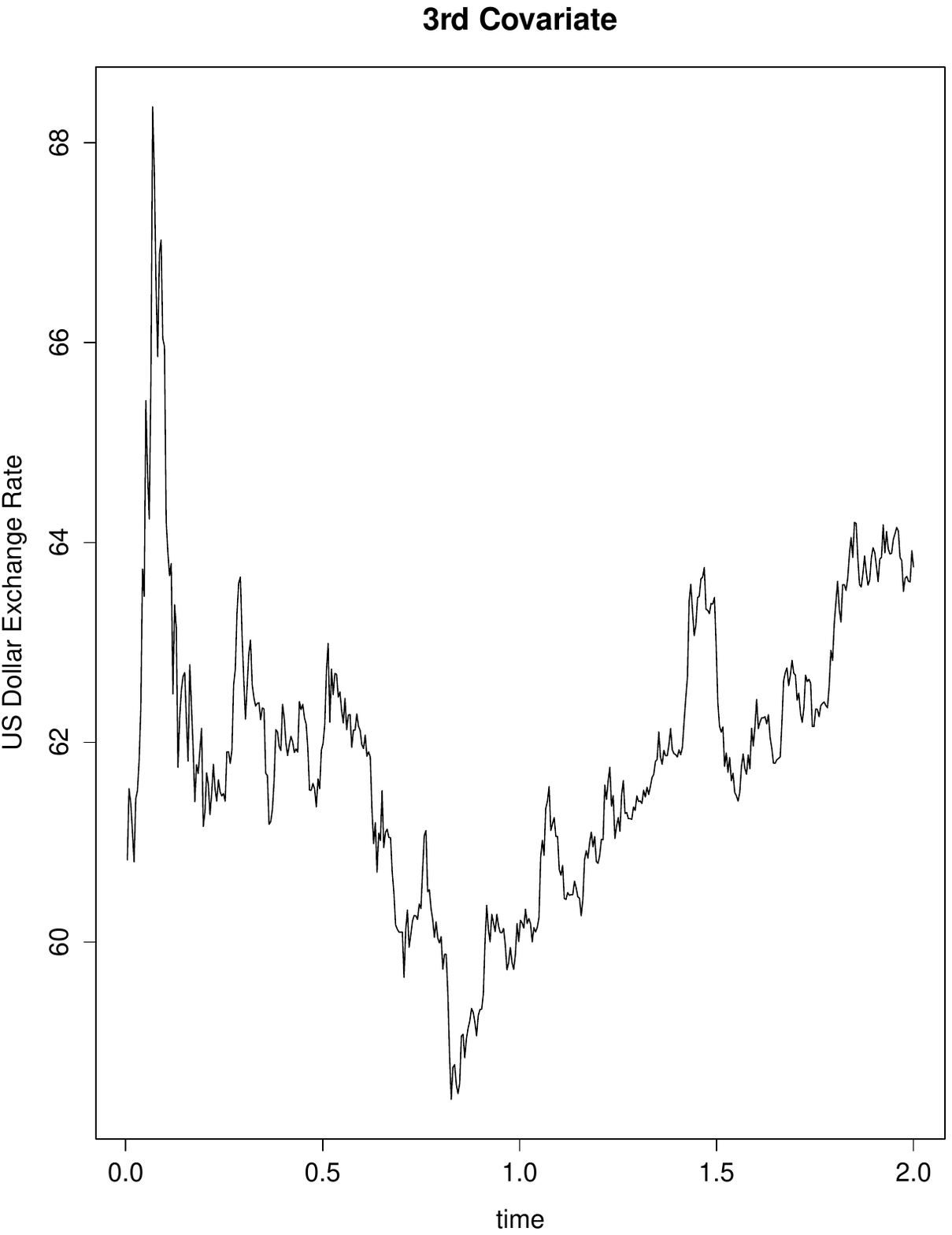}
\label{fig:realdata_sfig4}
\end{subfigure}
\caption{Covariates.}
\label{fig:cov}
\end{figure}

Hence, we collect the values of the aforementioned time varying covariates 
during the time range August $5$, 2013, to June $30$, 2015. The pattern of the covariates are 
displayed in Figure \ref{fig:cov}.

We denote these three covariates by $c_1,c_2,c_3$ respectively. Now, our considered $SDE$ models
for national stock exchange data associated with the $15$ companies are the
following:
\begin{equation}
dX_i(t)=(\theta_1^i+\theta_2^ic_1(t)+\theta_3^ic_2(t)+\theta_4^ic_3(t))(\theta_5^i+\theta_6^iX_i(t))dt+A^iX_i(t)^{B^i}dW_i(t),
\label{eq:sde4_appl}
\end{equation}
for $i=1,\cdots,15$. 

\subsection{Selection of covariates by Bayes factor}
\label{subsec:covariate_bf}

Among the considered three time varying covariates we now select the best set of covariate combinations
for the $15$ companies among $100$ such sets through Bayes factor, 
computing the log-marginal-likelihoods with respect to the normal prior
on the parameter set, assuming {\it a priori} independence of the parameter components 
with individual means being the corresponding $MLE$ (based on simulated annealing)
and $0.01^2$ variance (relatively small variance ensured numerical stability of th marginal likelihood).
Table \ref{table:companies} provides the sets of covariates for the $15$ companies obtained
by our Bayes factor analysis. Also observe that each of the three covariates occurs about $50\%$ times
among the companies, demonstrating that overall impact of these on national stock exchange is
undeniable.
\begin{table}[h]
\centering
\caption{Company-wise covariates obtained by Bayes factor analysis}
\label{table:companies}
\begin{tabular}{|c||c|}
\hline
Company & Covariates\\
\hline
$1$ & Bank rate\\ 
$2$ & US dollar exchange rate\\
$3$ & None\\
$4$ & None\\
$5$ & Bank rate and US dollar exchange rate\\
$6$ & Bank rate and US dollar exchange rate\\
$7$ & IIP general index and US dollar exchange rate\\
$8$ & Bank rate\\
$9$ & IIP general index and Bank rate\\
$10$ & IIP general index\\
$11$ & IIP general index, Bank rate and US dollar exchange rate\\
$12$ & IIP general index and Bank rate\\
$13$ & US dollar exchange rate\\
$14$ & IIP general index, Bank rate and US dollar exchange rate\\
$15$ & IIP general index\\
\hline
\end{tabular}
\end{table}

\section{Summary and discussion}
\label{sec:conclusion}

This article establishes the asymptotic theory of Bayes factors when the models to be compared are
systems of $SDE$'s consisting of time-dependent covariates and random effects, assuming that the
number of individuals as well as the domains of observations of the individuals increase
indefinitely. Different initial values for different $SDE$'s are also allowed. 
The only instance of related effort in this direction is that of \ctn{Maitra15a}. 
The main difference of our undertaking with that of \ctn{Maitra15a} is that
they assumed the domains of observations to be fixed for the individuals, a consequence
being that incorporation of random effects in their model was not possible from the asymptotic
perspective. Moreover, in their case, a single set of covariates was associated with all the individuals,
but here our random effects set-up allows different sets of time-dependent covariates for different individuals.

To proceed, we first needed to build
an asymptotic theory of Bayes factors for comparing two individual $SDE$'s, rather than two systems of $SDE$'s,
as the domain of observation expands. Our results in this regard, which help
formulate our asymptotic theory for comparing two systems of $SDE$'s using Bayes factors, 
are perhaps also of independent interest, being possibly the first ever results in this direction of research.
Although the relevant variance did not converge to zero when two individual $SDE$'s are compared,
we are able to establish almost sure exponential convergence of the Bayes factor when the number of subjects are allowed to
increase indefinitely.
Importantly, our theory covers both $iid$ and non-$iid$ cases.

Our simulation studies associated with covariate selection demonstrate that Bayes factor yields
consistent results even in non-asymptotic situations. Bayes factor analysis of a real data on company-wise national
stock exchange also yielded plausible sets of covariates for the companies.


Note that our current asymptotic Bayes factor theory remains valid
for comparison between $iid$ and non-$iid$ models.
For instance, if the true model consititutes an $iid$ system, then $f_{0i}\equiv f_0\equiv f_{\btheta_0}$; 
the rest remains the same as the theory for our non-$iid$ setting. The situation is analogous when the other model forms
an $iid$ system.

\section*{Acknowledgments}
We are thankful to Dr. Diganta Mukherjee for drawing our attention towards the website 
{\it www.nseindia.com} at which the real data was available. 
The first author gratefully acknowledges her CSIR Fellowship, Govt. of India.

\newpage
\input{supp_T}

\normalsize
\bibliographystyle{natbib}
\bibliography{irmcmc}

\end{document}

%% file: supp_T.tex
\renewcommand\thefigure{S-\arabic{figure}}
\renewcommand\thetable{S-\arabic{table}}
\renewcommand\thesection{S-\arabic{section}}

\setcounter{section}{0}
\setcounter{figure}{0}
\setcounter{table}{0}

\begin{center}
{\bf \Large Supplementary Material}
\end{center}

Throughout, we refer to our main manuscript as MB.

\section{Proof of Lemma 2 of MB}
\label{sec:proof_lemma_14}
Due to compactness of $\bGamma$ it follows, using the form of $\phi_{\bxi_j}$ provided in
(H5$^{\prime\prime}$), that the convergences (3.9), (3.10) and (3.11) of MB 
are uniform over $\bGamma$. The same form shows that the above integrals  
are continuous in $\bxi_1$, for every $T>0$. Hence, due to uniform convergence, the limits 
$\bar\phi^{(1)}_{\bxi_1}$, $\bar\phi^{(2)}_{\bxi_1}$ and $\bar\phi^{(2)}_{\bxi_0,\bxi_1}$
are continuous in $\bxi_1$.

\section{Proof of Lemma 3 of MB}
\label{sec:proof_lemma_15}
The proofs of (i) -- (iv) follow from (H1$^{\prime\prime}$), (H4$^{\prime\prime}$), the results
(3.10) and (3.11) of MB
following from (H5$^{\prime\prime}$), 
(3.1) and its asymptotic form (3.2) 
(with $k=2$), 
using the relation (2.8) of MB. 
To prove (v), 
note that, since for any $k\geq 1$, 
it holds, due to (H4$^{\prime\prime}$), (3.1) of MB 
and boundedness of $\phi_{\bxi_j}$ on $[0,T]$, that for $j=0,1$,
$$E\int_0^T\left|\frac{\phi_{\bxi_j}(s)b_{\bbeta_j}(s,X(s))}{\sigma(s,X(s))}\right|^{2k}ds<\infty,$$
it follows from Theorem 7.1 of \ctn{Mao11}, page 39, that
\begin{align}
E\left|\int_0^T\frac{\phi_{\bxi_j}(s)b_{\bbeta_j}(s,X(s))}{\sigma(s,X(s))}dW(s)\right|^{2k}
\leq\left(k(2k-1)\right)^{k}T^{k-1}E\int_0^T
\left|\frac{\phi_{\bxi_j}(s)b_{\bbeta_j}(s,X(s))}{\sigma(s,X(s))}\right|^{2k}ds.
\end{align}
Hence, using Chebychev's inequality, it follows that for any $\epsilon>0$,
\begin{align}
&P\left(\left|\frac{1}{T}\int_0^T\frac{\phi_{\bxi_j}(s)b_{\bbeta_j}(s,X(s))}{\sigma(s,X(s))}dW(s)\right|
>\epsilon\right)\notag\\
&\quad\quad\quad<\epsilon^{-2k}\left(k(2k-1)\right)^{k}T^{-(k+1)}E\int_0^T
\left|\frac{\phi_{\bxi_j}(s)b_{\bbeta_j}(s,X(s))}{\sigma(s,X(s))}\right|^{2k}ds.
\end{align}
In particular, if $k=2$ is chosen, then it 
follows from the above inequality, (H4$^{\prime\prime}$) , (3.2) of MB, 
and boundedness of
$\phi_{\bxi_j}$ on $[0,T]$, that
\begin{equation*}
\sum_{T=1}^{\infty}P\left(\left|\frac{1}{T}
\int_0^T\frac{\phi_{\bxi_j}(s)b_{\bbeta_j}(s,X(s))}{\sigma(s,X(s))}dW(s)\right|>\epsilon\right)
<\infty,
\end{equation*}
proving that 
$$\frac{1}{T}\int_0^T
\frac{\phi_{\bxi_j}(s)b_{\bbeta_j}(s,X(s))}{\sigma(s,X(s))}dW(s)\stackrel{a.s.}{\longrightarrow} 0.$$

To prove (vi), first note that 
\begin{align}
&E\left|\frac{1}{T}\int_0^T\left[\frac{\phi^2_{\bxi_j}(s)b^2_{\bbeta_j}(s,X(s))}{\sigma^2(s,X(s))}
-\phi^2_{\bxi_j}(s)\kappa_j(\bbeta_j)\right]ds\right|^{2k}\notag\\
&\leq T^{-1}E\int_0^T\phi^{4k}_{\bxi_j}(s)\left|\frac{b^2_{\bbeta_j}(s,X(s))}{\sigma^2(s,X(s))}
-\kappa_j(\bbeta_j)\right|^{2k}ds\notag\\
&\leq K_4T^{-1}E\int_0^T\left|\frac{b^2_{\bbeta_j}(s,X(s))}{\sigma^2(s,X(s))}
-\kappa_j(\bbeta_j)\right|^{2k}ds,
\end{align}
for some finite constant $K_4>0$.
The second last inequality is by H\"{o}lder's inequality, and the last inequality
holds because 
$\phi_{\bxi_j}(t)$ is uniformly bounded on $[0,\infty]$ thanks to compactness of $\boldcal Z$ and continuity of 
the functions $g_l;l=1,\ldots,p$. Hence, for any $\epsilon>0$,
\begin{align}
&P\left(\left|\frac{1}{T}\int_0^T\frac{\phi^2_{\bxi_j}(s)b^2_{\bbeta_j}(s,X(s))}{\sigma^2(s,X(s))}ds
-\kappa_j(\bbeta_j)\frac{1}{T}\int_0^T\phi^2_{\bxi_j}(s)ds\right|>\epsilon\right)\notag\\
&\qquad\qquad
<K_4\epsilon^{-2k}T^{-1}E\int_0^T\left|\frac{b^2_{\bbeta_j}(s,X(s))}{\sigma^2(s,X(s))}-\kappa_j(\bbeta_j)\right|^{2k}ds.
\notag
\end{align}
In the same way as the proof of (v), it follows, using the above inequality, 
(3.3) and (3.2) of MB, that
\begin{equation*}
\sum_{T=1}^{\infty}P\left(\left|\frac{1}{T}\int_0^T\frac{\phi^2_{\bxi_j}(s)b^2_{\bbeta_j}(s,X(s))}{\sigma^2(s,X(s))}ds
-\kappa_j(\bbeta_j)\frac{1}{T}\int_0^T\phi^2_{\bxi_j}(s)ds\right|>\epsilon\right)
<\infty.
\end{equation*}
That is, 
$$\frac{1}{T}\int_0^T\frac{\phi^2_{\bxi_j}(s)b^2_{\bbeta_j}(s,X(s))}{\sigma^2(s,X(s))}ds
-\kappa_j(\bbeta_j)\frac{1}{T}\int_0^T\phi^2_{\bxi_j}(s)ds\rightarrow 0,$$
almost surely, as $T\rightarrow\infty$. Since, as $T\rightarrow\infty$, 
$\frac{1}{T}\int_0^T\phi^2_{\bxi_j}(s)ds\rightarrow\bar\phi^{(2)}_{\bxi_j}$
by (3.10) of MB, 
the result follows.
Using (3.4) instead of (3.3), 
(vii) can be proved in the same way as (vi). 
The proofs of 
(viii) and (ix) follow from (v), (vi) and (vii), using the relation (2.8). 

\section{Proof of Lemma 4 of MB}
\label{sec:proof_lemma_16}
Using the Cauchy-Schwartz inequality twice we obtain
\begin{align}
&\frac{1}{T}\int_0^T
E\left(\frac{\phi_{\bxi_0}(s)\phi_{\bxi_1}(s)b_{\bbeta_0}(s,X(s))b_{\bbeta_1}(s,X(s))}{\sigma^2(s,X(s))}\right)ds
\notag\\
&\leq\frac{1}{T}\int_0^T
\sqrt{E\left(\frac{\phi^2_{\bxi_0}(s)b^2_{\bbeta_0}(s,X(s))}{\sigma^2(s,X(s))}\right)}
\times \sqrt{E\left(\frac{\phi^2_{\bxi_1}(s)b^2_{\bbeta_1}(s,X(s))}{\sigma^2(s,X(s))}\right)}ds\notag\\
&=\sqrt{\kappa_0(\bbeta_0)\kappa_1(\bbeta_1)+O\left(\exp\left(T^3-T^5\right)\right)}\times
\frac{1}{T}\int_0^T\left|\phi_{\bxi_0}(s)\right|\times\left|\phi_{\bxi_1}(s)\right|ds\notag\\
&\leq \sqrt{\kappa_0(\bbeta_0)\kappa_1(\bbeta_1)+O\left(\exp\left(T^3-T^5\right)\right)}\times
\left(\sqrt{\frac{1}{T}\int_0^T\phi^2_{\bxi_0}(s)ds}\right)\times
\left(\sqrt{\frac{1}{T}\int_0^T\phi^2_{\bxi_1}(s)ds}\right).
\label{eq:phi_kappa_2}
\end{align}
Taking the limit of both sides of (\ref{eq:phi_kappa_2}) as $T\rightarrow\infty$, using (ii) of 
Lemma 3 
and the limits (3.10), 
the result 
follows.

\section{Proof of Lemma 5 of MB}
\label{sec:proof_lemma_17}
For any $h\in (0,t)$,
\begin{align}
\frac{I_t}{I_{t-h}} 
&=\frac{\int_{\bTheta} \frac{f_{\btheta_1,t}(X_{0,t})}{f_{\btheta_0,t}(X_{0,t})}\pi (d\btheta_1)}
{\int_{\bTheta} \frac{f_{\btheta_1,{t-h}}(X_{0,t-h})}{f_{\btheta_0,{t-h}}(X_{0,t-h})}\pi (d\btheta_1)}\notag\\
&=\frac{\int_{\bTheta} \exp\left (\left(U_{\btheta_1,t} 
- U_{\btheta_0,t}\right)-\frac{\left(V_{\btheta_1,t}
-V_{\btheta_0,t}\right)}{2}\right) \pi(d\btheta_1)}
{\int_{\bTheta} \exp \left(\left(U_{\btheta_1,t-h} 
- U_{\btheta_0,t-h}\right)-\frac{\left(V_{\btheta_1,t-h}
-V_{\btheta_0,t-h}\right)}{2}\right) \pi(d\btheta_1)}\notag\\
&=\int_{\bTheta} \exp\left(\left(U_{\btheta_1,t-h} 
- U_{\btheta_0,t-h}\right)-\frac{\left(V_{\btheta_1,t-h}
-V_{\btheta_0,t-h}\right)}{2}\right)\notag\\
&\quad\times\frac{\exp\left(U_{\btheta_1,t-h,t}-U_{\btheta_0,t-h,t}
-\frac{\left(V_{\btheta_1,t-h,t}
-V_{\btheta_0,t-h,t}\right)}{2}\right)\pi(d\btheta_1)}
{\int_{\bTheta}\exp\left(\left(U_{\btheta_1,t-h} - U_{\btheta_0,t-h}\right)
-\frac{(V_{\btheta_1,t-h}-V_{\btheta_0,t-h})}{2}\right)\pi(d\btheta_1)}\notag\\
&=E\left[\exp\left(\int_{t-h}^t\frac{\phi_{\bxi_1}(s)b_{\bbeta_1}(s,X(s))}{\sigma^2(s,X(s))}dX(s)
-\int_{t-h}^t\frac{\phi_{\bxi_0}(s)b_{\bbeta_0}(s,X(s))}{\sigma^2(s,X(s))}dX(s)\right.\right.\notag\\
&\left.\left.\quad-\left(\frac{1}{2}
\int_{t-h}^t\frac{\phi_{\bxi_1}^2(s)b_{\bbeta_1}^2(s,X(s))}{\sigma^2(s,X(s))}ds
-\frac{1}{2}\int_{t-h}^t\frac{\phi_{\bxi_0}^2(s)b_{\bbeta_0}^2(s,X(s))}{\sigma^2(s,X(s))}ds\right)\right)
\bigg\vert\mathcal F_{t-h}\right]\notag\\
&=\frac{E\left[\exp\left(\int_{t-h}^t\frac{\phi_{\bxi_1}(s)b_{\bbeta_1}(s,X(s))}{\sigma^2(s,X(s))}dX(s)
-\frac{1}{2}\int_{t-h}^t\frac{\phi_{\bxi_1}^2(s)b_{\bbeta_1}^2(s,X(s))}{\sigma^2(s,X(s))}ds\right)
\bigg |\mathcal F_{t-h}\right]}{\exp\left(\int_{t-h}^t\frac{\phi_{\bxi_0}(s)b_{\bbeta_0}(s,X(s))}{\sigma^2(s,X(s))}dX(s)
-\frac{1}{2}\int_{t-h}^t\frac{\phi_{\bxi_0}^2(s)b_{\bbeta_0}^2(s,X(s))}{\sigma^2(s,X(s))}ds\right)}\notag\\
&=\frac{\hat f_{t-h,t}(X_{t-h,t})}{f_{\btheta_0,t-h,t}(X_{t-h,t})}.
\end{align}
Hence, the result 
holds.

\section{Proof of Theorem 6 of MB}
\label{sec:proof_theorem_18}
Let us consider
\begin{equation}
S_{Tq_n}=Tq_n\sum_{r=0}^{n(T)-1}\left(J'_{rTq_n}+\tilde{\mathcal K}'_{rTq_n}\right), 
\label{eq:Rsum}
\end{equation}
where $q_n=\frac{1}{n(T)}$, where, given $T>0$, $n(T)$ is the 
number of intervals partitioning $[0,T]$ each of length $\frac{T}{n(T)}$. We assume that as $T\rightarrow\infty$,
$\frac{T}{n(T)}\rightarrow 0$.

It follows, using (H9$^{\prime\prime}$), 
that for any $T>0$, 
\begin{align}
E\left(\frac{S_{Tq_n}}{T}\right)&\rightarrow \frac{1}{T}\int_0^T\frac{d}{dt}E_{\btheta_0}\left(\log I_t\right)dt
+\frac{1}{T}\int_0^TE_{\btheta_0}\left[\tilde{\mathcal K}'_t(f_{\btheta_0},\hat f)\right]dt\notag\\
&=E_{\btheta_0}\left(\frac{1}{T}\log I_T\right)
+\frac{1}{T}\int_0^TE_{\btheta_0}\left[\tilde{\mathcal K}'_t(f_{\btheta_0},\hat f)\right]dt,
\label{eq:conv1}
\end{align}
as $n(T)\rightarrow\infty$, for any given $T>0$.
Also, since due to (4.26) of MB, 
$E\left(J'_{rTq_n}+\tilde{\mathcal K}'_{rTq_n}|\mathcal F_{rTq_n}\right)=0$,
we must have $E\left(J'_{rTq_n}+\tilde{\mathcal K}'_{rTq_n}\right)
=E\left[E\left(J'_{rTq_n}+\tilde{\mathcal K}'_{rTq_n}|\mathcal F_{rTq_n}\right)\right]=0$, for any $r,T,n(T)$.
Hence, $E\left(\frac{S_{Tq_n}}{T}\right)=0$ for any $T,n(T)$. 
Thus, it follows from (\ref{eq:conv1}), that
\begin{equation}
E_{\btheta_0}\left(\frac{1}{T}\log I_T\right)
+\frac{1}{T}\int_0^TE_{\btheta_0}\left[\tilde{\mathcal K}'_t(f_{\btheta_0},\hat f)\right]dt
\rightarrow 0,\quad\mbox{as}~T\rightarrow\infty.
\label{eq:conv2}
\end{equation}
We now deal with the second term of the left hand side of (\ref{eq:conv2}).
Since, by (H6$^{\prime\prime}$), 
$$\pi\left(\btheta_1: \underset{t}{\inf}~\tilde{\mathcal K}'_t(f_{\btheta_0},f_{\btheta_1})\geq \delta\right)=1,$$
it holds that $\tilde{\mathcal K}'_t(f_{\btheta_0},f_{\btheta_1})\geq \delta$ for all $t$ with probability 1,
so that 
\[
\tilde{\mathcal K}'_t(f_{\btheta_0},\hat f)=\tilde{\mathcal K}'_t(f_{\btheta_0},\hat f_{A_t(\delta)}),
\]
where $A_t(\delta)$ is given by (4.19) of MB. 
The $Q^*$ property 
implies that  
\begin{equation}
\underset{T}{\lim\inf}~ \frac{1}{T}\int_0^TE_{\btheta_0}\left[\tilde{\mathcal K}'_t(f_{\btheta_0},\hat f)\right]dt\geq \delta.
\label{eq:liminf_K_average_T}
\end{equation}
The results (\ref{eq:conv2}) and (\ref{eq:liminf_K_average_T}) imply that
\begin{align}
\underset{T}{\limsup}~ E_{\btheta_0}\left(\frac{1}{T}\log I_T\right)\leq -\delta.
\label{eq:limsup_T}
\end{align}
Now observe that
\begin{align}
I_T&=\int_{\bTheta} \exp\left(U_{\btheta_1,T}
-U_{\btheta_0,T}\right)
\times\exp\left\{-\frac{1}{2}\left(V_{\btheta_1,T}
-V_{\btheta_0,T}\right)\right\}
\pi(d\btheta_1)\notag\\
&\geq 
\int_{\mathcal N_0(c)}\exp\left\{T\left(\frac{U_{\btheta_1,T}}{T}
-\frac{U_{\btheta_0,T}}{T}\right)\right\}\notag\\
&\quad\quad\times\exp\left\{-\frac{T}{2}\left(\frac{V_{\btheta_1,T}}{T}
-\frac{V_{\btheta_0,T}}{T}\right)\right\}
\pi(d\btheta_1),
\label{eq:lower_bound1}
\end{align}
where $c>0$, and
\begin{align}
\mathcal N_0(c)&=\left\{\btheta_1\in\bTheta:\delta\leq\bar{K}^{\infty}\left(f_{\btheta_0},f_{\btheta_1}\right)
\leq\delta+c\right\}\notag\\
&=\left\{\btheta_1\in\bTheta:\delta\leq\frac{\bar\phi^{(2)}_{\bxi_0}}{2}\kappa_0(\bbeta_0)
-\bar\phi^{(2)}_{\bxi_0,\bxi_1}\bar{\kappa}(\bbeta_0,\bbeta_1)
+\frac{\bar\phi^{(2)}_{\bxi_1}}{2}\kappa_1(\bbeta_1)
\leq\delta+c\right\},\notag
\end{align}
the second line following from (4.15) of MB. 
Using Jensen's inequality, we obtain
\begin{align}
\frac{1}{T}\log\left(I_T\right)&\geq 
\int_{\mathcal N_0(c)}\left[\left(\frac{U_{\btheta_1,T}}{T}
-\frac{U_{\btheta_0,T}}{T}\right)\right.\notag\\
&\quad\quad\left.-\frac{1}{2}\left(\frac{V_{\btheta_1,T}}{T}
-\frac{V_{\btheta_0,T}}{T}\right)\right]
\pi(d\btheta_1).
\label{eq:lower_bound2}
\end{align}
By (vi) -- (ix) of Lemma 3 of MB, 
the integrand of the right hand side of the above
inequality, which we denote by $g_{X_T}(\btheta_1)$, converges to 
$g(\btheta_1)=-\left[\frac{\bar\phi^{(2)}_{\bxi_1}}{2}\kappa_1(\bbeta_1)
-\bar\phi^{(2)}_{\bxi_0,\bxi_1}\bar{\kappa}(\bbeta_0,\bbeta_1)
+\frac{\bar\phi^{(2)}_{\bxi_0}}{2}\kappa_0(\bbeta_0)\right]$, pointwise for every $\btheta_1$, given any 
path of the process $X$ in the complement of the null set. 
Due to (H1$^{\prime\prime}$), (H4$^{\prime\prime}$) and (3.2) of MB, 
$\underset{T}{\sup}~E_{\btheta_1}\left[g_{X_T}(\btheta_1)\right]^2<\infty$,
so that $\left\{g_{X_T}(\btheta_1)\right\}_{T=1}^{\infty}$ is uniformly integrable. Hence,
$$\int_{\mathcal N_0(c)}g_{X_T}(\btheta_1)\pi(d\btheta_1)\rightarrow 
\int_{\mathcal N_0(c)}g(\btheta_1)\pi(d\btheta_1),$$
given any path of the process $X$ in the complement of the null set.
Let us denote the left hand side of the above by $H_{X_T}$ let $H$ denote the right hand side.
We just proved that $H_{X_T}$ converges to $H$ almost surely. Now observe that
\begin{align}
\underset{T}{\sup}~E_{\btheta_0}\left[H_{X_T}\right]^2
&=\underset{T}{\sup}~E_{\btheta_0}\left[\int_{\mathcal N_0(c)}g_{X_T}(\btheta_1)\pi(d\btheta_1)\right]^2\notag\\
&\leq \int_{\mathcal N_0(c)}\underset{T}{\sup}~E_{\btheta_0}\left[g_{X_T}(\btheta_1)\right]^2\pi(d\btheta_1).\notag
\end{align}
Again, due to (H1$^{\prime\prime}$), (H4$^{\prime\prime}$) and (3.2) of MB, 
the last expression is finite, proving uniform integrability
of $\left\{H_{X_T}\right\}_{T=1}^{\infty}$. Hence, 
$$\underset{T\rightarrow\infty}{\lim}~E_{\btheta_0}\left(H_{X_T}\right)= E_{\btheta_0}\left(H\right)=H.$$
It follows that
\begin{align}
&\underset{T\rightarrow\infty}{\liminf}~E_{\btheta_0}\left[\frac{1}{T}\log\left(I_T\right)\right]
\geq \underset{T\rightarrow\infty}{\liminf}~E_{\btheta_0}\left(H_{X_T}\right)
=E_{\btheta_0}\left(H\right)\notag\\
&=\int_{\mathcal N_0(c)}g(\btheta_1)\pi(\btheta_1)d\btheta_1\notag\\
&\geq-\left(\delta+c\right)\pi\left(\mathcal N_0(c)\right)\notag\\
&\geq-\left(\delta+c\right).
\end{align}
Since the above holds for arbitrary $c>0$,
it holds that
\begin{equation}
\underset{T\rightarrow\infty}{\lim\inf}~E_{\btheta_0}\left(\frac{1}{T}\log I_T\right)\geq -\delta.
\label{eq:liminf_T}
\end{equation}
Thus (\ref{eq:limsup_T}) and (\ref{eq:liminf_T}) together help us conclude that 
$$E_{\btheta_0}\left(\frac{1}{T}\log I_T\right)\rightarrow-\delta,$$ as $T\rightarrow\infty$.

We now show that the variance of $\frac{1}{T}\log I_T$ is $O(1)$, as $T\rightarrow\infty$.
First note, due to compactness of $\bTheta$, the mean value theorem for integrals
ensure existence of $\acute\btheta_1=(\acute\bbeta_1,\acute\bxi_1)=\acute\btheta_1(W)\in\bTheta$, depending 
on the Wiener process $W$ 
such that
\begin{equation}
\log I_T=\left(U_{\acute\btheta_1,T}
-U_{\btheta_0,T}\right)
-\frac{1}{2}\left(V_{\acute\btheta_1,T}
-V_{\btheta_0,T}\right).
\label{eq:log_I_T}
\end{equation}

Now note that the results presented in Lemma 3 of MB continue to hold even when
$\btheta_1$ is replaced with $\acute\btheta_1$. Specifically, the following hold
in addition to the results of Lemma 3:
\begin{align}
\quad &\frac{V_{\acute\btheta_1,T}}{T}-\bar \phi^{(2)}_{\acute\bxi_1}\kappa_1(\acute\bbeta_1)
\stackrel{a.s.}{\longrightarrow} 0; \notag\\ 
\quad &\frac{U_{\acute\btheta_1,T}}{T}-\bar\phi^{(2)}_{\bxi_0,\acute\bxi_1}\bar{\kappa}(\bbeta_0,\acute\bbeta_1)
\stackrel{a.s.}{\longrightarrow}0.\notag 
\end{align}
It follows that, as $T\rightarrow\infty$,
$$\ell_T=\frac{1}{T}\log I_T-\bar\phi^{(2)}_{\bxi_0,\acute\bxi_1}\bar{\kappa}(\bbeta_0,\acute\bbeta_1)
+\frac{1}{2}\bar \phi^{(2)}_{\acute\bxi_1}\kappa_1(\acute\bbeta_1)
+\frac{1}{2}\bar \phi^{(2)}_{\bxi_0}\kappa_0(\bbeta_0)
\stackrel{a.s.}{\longrightarrow} 0.$$
By uniform integrability arguments, which follow in similar lines as the proofs of Lemmas 1 and 10 of \ctn{Maitra15a} using 
(H4$^{\prime\prime}$), compactness, and Cauchy-Schwartz, it holds that
\begin{equation}
Var_{\btheta_0}\left(\ell_T\right)\rightarrow 0,~\mbox{as}~T\rightarrow\infty.
\label{eq:var_ell}
\end{equation}
Hence,
\begin{align}
Var_{\btheta_0}\left(\frac{1}{T}\log I_T\right)
&=Var_{\btheta_0}\left(\ell_T+\bar\phi^{(2)}_{\bxi_0,\acute\bxi_1}\bar{\kappa}(\bbeta_0,\acute\bbeta_1)
-\frac{1}{2}\bar \phi^{(2)}_{\acute\bxi_1}\kappa_1(\acute\bbeta_1)
-\frac{1}{2}\bar \phi^{(2)}_{\bxi_0}\kappa_0(\bbeta_0)\right)\notag\\
&= Var_{\btheta_0}\left(\ell_T\right)
+Var_{\btheta_0}\left(\bar\phi^{(2)}_{\bxi_0,\acute\bxi_1}\bar{\kappa}(\bbeta_0,\acute\bbeta_1)
-\frac{1}{2}\bar \phi^{(2)}_{\acute\bxi_1}\kappa_1(\acute\bbeta_1)
-\frac{1}{2}\bar \phi^{(2)}_{\bxi_0}\kappa_0(\bbeta_0)\right)\notag\\
&\qquad+2Cov_{\btheta_0}\left(\ell_T,\bar\phi^{(2)}_{\bxi_0,\acute\bxi_1}\bar{\kappa}(\bbeta_0,\acute\bbeta_1)
-\frac{1}{2}\bar \phi^{(2)}_{\acute\bxi_1}\kappa_1(\acute\bbeta_1)
-\frac{1}{2}\bar \phi^{(2)}_{\bxi_0}\kappa_0(\bbeta_0)\right).
\label{eq:var2}
\end{align}
By (\ref{eq:var_ell}), the first term of (\ref{eq:var2}) goes to zero as $T\rightarrow\infty$, and
the third, covariance term tends to zero by Cauchy-Schwartz and (\ref{eq:var_ell}). In other words,
as $T\rightarrow\infty$,
\begin{equation}
\left|Var_{\btheta_0}\left(\frac{1}{T}\log I_T\right)
-Var_{\btheta_0}\left(\bar\phi^{(2)}_{\bxi_0,\acute\bxi_1}\bar{\kappa}(\bbeta_0,\acute\bbeta_1)
-\frac{1}{2}\bar \phi^{(2)}_{\acute\bxi_1}\kappa_1(\acute\bbeta_1)\right)\right|
\rightarrow 0.
\label{eq:var3}
\end{equation}
However, $$Var_{\btheta_0}\left(\bar\phi^{(2)}_{\bxi_0,\acute\bxi_1}\bar{\kappa}(\bbeta_0,\acute\bbeta_1)
-\frac{1}{2}\bar \phi^{(2)}_{\acute\bxi_1}\kappa_1(\acute\bbeta_1)\right)\nrightarrow 0,$$
unless $\bar\phi^{(2)}_{\bxi_0,\acute\bxi_1}\bar{\kappa}(\bbeta_0,\acute\bbeta_1)
-\frac{1}{2}\bar \phi^{(2)}_{\acute\bxi_1}\kappa_1(\acute\bbeta_1)$ is constant almost surely.
It then follows from (\ref{eq:var3}) that
\begin{equation}
Var_{\btheta_0}\left(\frac{1}{T}\log I_T\right)=O(1),~\mbox{as}~T\rightarrow\infty.
\label{eq:var4}
\end{equation}

\section{Proof of Theorem 8 of MB}
\label{sec:proof_theorem_20}
In our set-up it follows from (2.7) of MB that 
\begin{equation}
\frac{1}{nT}\log\tilde I_{n,T}
=\frac{1}{n}\sum_{i=1}^n\frac{1}{T}\log I_{i,T}.
\label{eq:bf_convergence_iid}
\end{equation}
In the $iid$ case, given $T>0$, 
using the above form, it follows by the strong law of large numbers, that 
\begin{equation}
\underset{n\rightarrow\infty}{\lim}\frac{1}{nT}\log\tilde I_{n,T}= E_{\btheta_0}\left(\frac{1}{T}\log I_{i,T}\right),
\label{eq:convergence_given_T}
\end{equation}
almost surely. Now, in the $iid$ situation, for each $i$,
$E_{\btheta_0}\left(\frac{1}{T}\log I_{i,T}\right)\rightarrow -\delta$, as $T\rightarrow\infty$.
Hence, taking limit as $T\rightarrow\infty$ on both sides of (\ref{eq:convergence_given_T}) yields
\begin{equation}
\underset{T\rightarrow\infty}{\lim}\underset{n\rightarrow\infty}{\lim}\frac{1}{nT}\log\tilde I_{n,T}
= \underset{T\rightarrow\infty}{\lim}E_{\btheta_0}\left(\frac{1}{T}\log I_{i,T}\right)=-\delta,
\label{eq:convergence_T_n}
\end{equation}
almost surely, proving the theorem.

\section{Proof of Lemma 10 of MB}
\label{sec:proof_lemma_22}
Note that, due to compactness of $\mathfrak X$ and $\boldcal Z$ and continuity of the covariates
in time $t$, there exists $x^*\in\mathfrak X$ and $\bz^*\in\mathfrak Z$, such that
\begin{equation}
\underset{x\in\mathfrak X,\bz\in\mathfrak Z}{\sup}~
\left|\frac{1}{T}E_{\btheta_0}\left(\log I_{x,T,\bz}\right) +\delta(x,\bz)\right|
=\left|\frac{1}{T}E_{\btheta_0}\left(\log I_{x^*,T,\bz^*}\right) +\delta(x^*,\bz^*)\right|
\rightarrow 0,
\label{eq:uniform_convergence}
\end{equation}
as $T\rightarrow\infty$, where the convergence is due to (6.2) of MB. 
Also, $\frac{1}{T}E_{\btheta_0}\left(\log I_{x,T,\bz}\right)$ is clearly
continuous in $(x,\bz)$ for every $T>0$ (the proof of this follows in the same way as that of Theorem 5 of \ctn{Maitra14a}). 
Combining this with the uniform convergence (\ref{eq:uniform_convergence})
it follows that $\delta(x,\bz)$ is also continuous in $(x,\bz)$.

\section{Proof of Lemma 11 of MB}
\label{sec:proof_lemma_23}
Note that the limit (6.3) of MB 
can be represented as
\begin{equation}
\underset{n\rightarrow\infty}{\lim}~\frac{1}{n}\sum_{r=0}^{n-1}
\varrho_{\delta}\left(\frac{r}{n}\right),
\label{eq:Riemann_delta}
\end{equation}
where $\varrho_{\delta}:[0,1]\mapsto\mathbb R^+$ is some continuous function satisfying 
$\varrho_{\delta}\left(\frac{r}{n}\right)=\delta\left(x^{r+1},\bz_{r+1}\right)$ for $r=0,\ldots,n-1$.
For the remaining points $y\in[0,1]$, we set $\varrho_{\delta}(y)=\delta(x,\bz)$,
where $(x,\bz)\in\mathfrak X\times\mathfrak Z$ is such that $\varrho_{\delta}(y)$ is continuous in $y\in[0,1]$. Since 
$\delta(x,\bz)$ is continuous in $(x,\bz)$, $\varrho_{\delta}(y)$ can be thus constructed.
Note that, it is possible to relate $y\in[0,1]$ to $(x,\bz)\in\mathfrak X\times\mathfrak Z$ by some continuous mapping
$G:\mathfrak X\times\mathfrak Z\mapsto[0,1]$, taking $(x,\bz)$ to $y$.
Thus, $\delta^{\infty}$ in (6.3) of MB
is the limit of the Riemann sum (\ref{eq:Riemann_delta}) associated with the continuous
function $\varrho_{\delta}$; the limit is given by the integral $\int_0^1\varrho_{\delta}(y)dy$. 
Since the domain of integration is $[0,1]$, it follows, using continuity of $\varrho_{\delta}$, 
that the integral is finite. Observe that for any given sequence $\left\{(x^i,\bz_i)\right\}_{i=1}^{\infty}$, one
can construct a continuous function $\varrho_{\delta}$ such that $\delta^{\infty}=\int_0^1\varrho_{\delta}(y)dy$.
In other words,
$\delta^{\infty}$ 
exists for all sequences $\left\{(x^i,\bz_i)\right\}_{i=1}^{\infty}$.

\section{Proof of Theorem 12 of MB}
\label{sec:proof_theorem_24}
For given $T>0$, it follows from (\ref{eq:var4}), compactness of $\bTheta$, $\mathfrak X$, $\boldcal Z$,
and continuity of the relevant functions $\phi_{\bxi_j}$, 
$b_{\bbeta_j}$, $g_1,\ldots,g_p$, $\kappa_0$, $\kappa_1$ and $\bar\kappa$, that
\begin{equation}
\underset{x\in\mathfrak X,\bz\in\mathfrak Z}{\sup}~Var_{\btheta_0}\left(\frac{1}{T}\log I_{x,T,\bz}\right)<\infty.
\label{eq:var_sup_Z}
\end{equation}
Hence, given $T>0$,
\begin{equation*}
\sum_{i=1}^{\infty}\frac{Var_{\btheta_0}\left(\frac{1}{T}\log I_{i,T}\right)}{i^2}<\infty.
\end{equation*}
It then follows due to Kolmogorov's strong law of large numbers for independent random variables that
\begin{equation}
\underset{n\rightarrow\infty}{\lim}~\frac{1}{n}\sum_{i=1}^n\left(\frac{1}{T}\log I_{i,T}\right)
=\underset{n\rightarrow\infty}{\lim}~\frac{1}{n}\sum_{i=1}^nE_{\btheta_0}\left(\frac{1}{T}\log I_{i,T}\right),
\label{eq:non_iid_convergence1}
\end{equation}
almost surely.

Now observe that the right hand side of (\ref{eq:non_iid_convergence1}) admits the following representation
\begin{equation}
\underset{n\rightarrow\infty}{\lim}~\frac{1}{n}\sum_{r=0}^{n-1}\breve\varrho\left(\frac{r}{n},T\right),
\label{eq:Riemann1}
\end{equation}
where $\breve\varrho(\cdot,T):[0,1]\mapsto\mathbb R$ is some continuous function depending upon $T$ with 
$\breve\varrho(\frac{r}{n},T)=E_{\btheta_0}\left(\frac{1}{T}\log I_{r+1,T}\right)$.
 
Since $E_{\btheta_0}\left(\frac{1}{T}\log I_{x,T,\bz}\right)$ is continuous in $(x,\bz)$, $\breve\varrho(y,T)$ can be constructed 
as in Lemma 11 of MB. 
Then, for almost all $y\in [0,1]$, $\breve\varrho(y,T)\rightarrow -\delta(x,\bz)$ as $T\rightarrow\infty$, for appropriate
$(x,\bz)$ associated with $y$ via $y=G(x,\bz)$ as in Lemma 11. 
Also, it follows from (\ref{eq:log_I_T}) that 
$\breve\varrho(\cdot,T)$ so constructed is uniformly bounded in $T>0$. 
Thus, the conditions
of the dominated convergence theorem are satisfied.

Since (\ref{eq:Riemann1}) is nothing but the Riemann sum associated with $\breve\varrho(\cdot)$, it follows that
\begin{equation}
\underset{n\rightarrow\infty}{\lim}~\frac{1}{n}\sum_{i=1}^nE_{\btheta_0}\left(\frac{1}{T}\log I_{i,T}\right)
=\int_0^1\breve\varrho(y,T)dy.
\label{eq:towards_DCT}
\end{equation}
By construction of $\breve\varrho(y,T)$, the dominated convergence theorem holds for the right hand side
of (\ref{eq:towards_DCT}). Hence, 
\begin{align}
&\underset{T\rightarrow\infty}{\lim}\underset{n\rightarrow\infty}{\lim}
~\frac{1}{n}\sum_{i=1}^nE_{\btheta_0}\left(\frac{1}{T}\log I_{i,T}\right)
=\underset{T\rightarrow\infty}{\lim}~\int_0^1\breve\varrho(y,T)dy
=\int_0^1\underset{T\rightarrow\infty}{\lim}~\breve\varrho(y,T)dy\notag\\
&\quad\quad=\underset{n\rightarrow\infty}{\lim}~\frac{1}{n}\sum_{i=1}^n
\underset{T\rightarrow\infty}{\lim}~E_{\btheta_0}\left(\frac{1}{T}\log I_{i,T}\right)
=-\underset{n\rightarrow\infty}{\lim}~\frac{1}{n}\sum_{i=1}^n\delta(x^i,\bz_i)=-\delta^{\infty}.
\label{eq:DCT}
\end{align}
Combining (\ref{eq:DCT}) with (\ref{eq:non_iid_convergence1}) it follows that
\begin{equation}
\underset{T\rightarrow\infty}{\lim}\underset{n\rightarrow\infty}{\lim}
~\frac{1}{n}\sum_{i=1}^n\left(\frac{1}{T}\log I_{i,T}\right)
=\underset{T\rightarrow\infty}{\lim}\underset{n\rightarrow\infty}{\lim}
~\frac{1}{n}\sum_{i=1}^nE_{\btheta_0}\left(\frac{1}{T}\log I_{i,T}\right)
=-\delta^{\infty},
\label{eq:non_iid_convergence2}
\end{equation}
almost surely.